\documentclass[preprint,12pt]{elsarticle}
\journal{Stochastics and Dynamics}

\usepackage {amssymb}
\usepackage {amsmath,amssymb,paralist,amsthm}
\usepackage {graphicx,color,epsfig,fancyhdr}
\usepackage {mathrsfs,bbm,dsfont}
\usepackage {hyperref,tikz}
\usepackage {latexsym,amsfonts}

\def\E{\mathbb{E}}
\def\R{\mathbb{R}}

\def\1{\mathbf{1}}

\def\P{\mathbb{P}}

\newtheorem{theorem}{Theorem}[section]
\newtheorem*{pf*}{Proof of Theorem 2.2}
\newtheorem{lemma}{Lemma}[section]

\begin{document}

\begin{frontmatter}



\title{Backward stochastic differential equations with unbounded generators}


\author{Bujar Gashi and Jiajie Li\footnote{Corresponding author.}}
\address{\scriptsize{Institute of Financial and Actuarial Mathematics (IFAM),\\
Department of Mathematical Sciences, The University of Liverpool, Liverpool, L69 7ZL, UK\\
Emails: Bujar.Gashi@liverpool.ac.uk; Jiajie.Li@me.com}}

\begin{abstract}
In this paper we consider two classes of backward stochastic differential equations. Firstly, under a Lipschitz-type condition on the generator of the equation, which can also be unbounded, we give sufficient conditions for the existence of a unique solution pair. The method of proof is that of Picard iterations and the resulting conditions are new. We also prove a comparison theorem. Secondly, under the linear growth and continuity assumptions on the possibly unbounded generator, we prove the existence of the solution pair. This class of equations is more general than
the existing ones.
\end{abstract}


\end{frontmatter}


\section{Introduction}
Let $(\Omega, \mathscr{F}, (\mathscr{F}_t,t\geq0),\mathds{P})$ be a given complete filtered probability space on which a $k$-dimensional standard
Brownian motion $(W(t), t\geq 0)$ is defined. We assume that $\mathscr{F}_t$ is the augmentation of $\sigma\{W(s):0\leq s\leq t\}$ by all the
$\mathds{P}$-null sets of $\mathscr{F}$. Consider the backward stochastic differential equation (BSDE):
\begin{equation}
\label{bsde}
y(t)= \xi+ \int_t^T f(s, y(s), z(s))ds-\int_t^T z(s)dW(s),\quad t\in[0,T],
\end{equation}
where $\xi$ is a given $\mathscr{F}_T$-measurable $\R^d$-valued random variable, and the generator $f: \Omega\times(0,T)\times\R^d\times\R^{d\times k}\rightarrow\R^d$ is a progressively measurable function.\\
Linear equations of the type (\ref{bsde}) were introduced by Bismut~\cite{bismut} in the context of stochastic linear quadratic control. The problem of existence and uniqueness of solution to the nonlinear equations (\ref{bsde}) was solved by Pardaoux and Peng~\cite{PP} under the global Lipschitz condition on $f$, i.e. under the assumption that there exists a real constant $c>0$ such that
\begin{eqnarray}
|f(t,y_1,z_1)-f(t,y_2,z_2)|\leq c(|y_1-y_2|+|z_1-z_2|),\label{lipschitz}
\end{eqnarray}
for all $y_1,y_2\in\R^d$, $z_1,z_2\in\R^{d\times k},$  $ (t,\omega)$  $a.e.$. Since then, the BSDEs have been studied extensively, and have found wide applicability in areas such as mathematical finance, stochastic control,  and stochastic controllability; see, for
example,~\cite{B},~\cite{NSM},~\cite{Ma},~\cite{M2},~\cite{XYZ},~\cite{peng},~\cite{G},~\cite{GP},~\cite{GP2}~\cite{wang}, and the references therein. One direction of
research has been to weaken the assumption of global Lipschitz condition (\ref{lipschitz}) by assuming only local Lipschitz condition (see~\cite{Ba}),
or non-Lipschitz condition of a particular form (see~\cite{M1},~\cite{WW}). In~\cite{KH} and ~\cite{DT} the authors permit for the generator function $f$ to be unbounded (they also consider a more general driving process than the Brownian motion). They assume the following global Lipschitz-type condition: there exist non-negative processes $c_1(\cdot)$ and $c_2(\cdot)$ such that
\begin{equation}
\label{lip}
|f(t,y_1,z_1)-f(t,y_2,z_2)|\leq c_1(t)|y_1-y_2|+c_2(t)|z_1-z_2|,
\end{equation}
for all $y_1, y_2\in\R^d$, $z_1, z_2\in\R^{d\times k}$, $ (t,\omega)$  $a.e.$. Clearly, this condition has great similarity with (\ref{lipschitz}).
However, different from (\ref{lipschitz}), here the processes $c_1(\cdot)$ and $c_2(\cdot)$ are not assumed to be bounded. The linear BSDEs with possibly unbounded coefficients are considered in~\cite{Yong1}, using a very different approach as compared to~\cite{KH},~\cite{DT}. The interest in these equations is not only theoretical, but is also motivated by applications in mathematical finance. Indeed, some very important interest rate models are given by stochastic differential equations (see, for
example,~\cite{Yong2},~\cite{DG},~\cite{BK}). The problem of market completeness (and thus of pricing and hedging of derivatives) in such models gives rise to BSDEs with possibly unbounded coefficients (see~\cite{Yong1} for details).

The first contribution of the present paper, as contained in the section \ref{LIP}, is to consider the problem of existence and uniqueness of a
solution pair $(y(\cdot),z(\cdot))$ for (\ref{bsde}) under condition (\ref{lip}). We do so under certain new conditions on the coefficients $c_1(\cdot)$, $c_2(\cdot)$, which are similar to those of \cite{KH}, but in general are not comparable. Moreover, our method of proof is different, since it is a modification of the Picard iteration procedure of~\cite{PP} rather than being based on a fixed point theorem as in \cite{KH}. We also give a comparison theorem for this class of equations. This generalises the classical result of Peng~\cite{Peng2},~\cite{Peng}, to the case of BSDEs with possibly unbounded coefficients.

Another important weakening of the assumptions on the generator, as compared to~\cite{PP},  was given in~\cite{LS} (see
also~\cite{HLP}). There it assumed that the generator is continuous with respect to $y$ and $z$, and it satisfies the linear growth condition
\begin{eqnarray}
|f(t,y,z)|\leq c(1+|y|+|z|),\label{ls1}
\end{eqnarray}
for all $y\in\R$, $z\in\R^{k}$, $ (t,\omega)$  $a.e.$. Under such conditions, it was shown that equation \eqref{bsde} admits a solution pair. In more recent papers~\cite{WW},~\cite{WH}, the linear growth condition \eqref{ls1} has been generalised to
\begin{eqnarray}
|f(t,y,z)|\leq c[q(t)+|y|+|z|],\label{ls2}
\end{eqnarray}
for all $y\in\R$, $z\in\R^{k}$, $ (t,\omega)$  $a.e.$. Here, different from \eqref{ls1}, the process $q(\cdot)$ is not assumed to be bounded.

The generator of BSDEs with a {\it quadratic} growth on the control variable $z$ satisfies the condition:
\begin{eqnarray}
|f(t,y,z)|\leq k_0+k_1|y|+k_2|z|^2.\label{q1}
\end{eqnarray}
An existence result for these equations was first given in~\cite{kob} where $k_0$ and $k_1$ are constant, $k_2$ is a given function of $y$, and the terminal value $\xi$ is assumed bounded. In~\cite{hu1},~\cite{hu2}, the assumption of a bounded $\xi$ was replaced with an integrability condition on the exponential of $\xi$, $k_0$ was permitted to be an unbounded process, whereas $k_1$ and $k_2$ were assumed constant (see further, for example,~\cite{BHM},~\cite{br2},~\cite{br1},~\cite{d2},~\cite{d3},~\cite{zheng}, where a nonlinear growth in $y$ is also permitted in some cases). In~\cite{k1},~\cite{k2},~\cite{mm},~\cite{marie}, the coefficient $k_1$ could also be an unbounded process. Under further assumptions on $f$, the uniqueness of the solution pair has been proved (see, for example,~\cite{br1},~\cite{d2},~\cite{d3}). We are not aware of an existence result for these types of BSDEs with coefficient $k_2$ being an unbounded process.

The second contribution of the present paper, which is contained in section \ref{CON}, is to consider a generator which is continuous in $y$ and $z$, but with a weaker linear growth condition than \eqref{ls2}. We assume that
\begin{eqnarray}
|f(t,y,z)|\leq c_0(t)+c_1(t)|y|+c_2(t)|z|,\label{ls3}
\end{eqnarray}
for all $y\in\R$, $z\in\R^{k}$, $ (t,\omega)$  $a.e.$. Here the processes $c_0(\cdot), c_1(\cdot), c_2(\cdot)$, are not assumed to be
bounded. By using the results of section 2, and appropriately modifying the approach of~\cite{BHM},~\cite{LS}, we prove the existence of a solution pair for \eqref{bsde}. Note that due to the unbounded nature of the process $c_2(t)$, conditions (\ref{q1}) and (\ref{ls3}) are in general not comparable if $k_2$ is a constant.

We conclude this introductory section with some notations:\\

$\bullet$ $|\cdot|$ is the Euclidian norm.\\

$\bullet$ $c_0(\cdot)$, $c_1(\cdot)$, $c_2(\cdot)$ are given $\R$-valued progressively measurable processes.\\

$\bullet$ $\gamma(\cdot)$, $\overline{\gamma(\cdot)}$ and $\widetilde{\gamma(\cdot)}$ are given $\R$-valued positive progressively measurable processes.\\

$\bullet$ $1<\beta_1, \widetilde{\beta_1}\in\R$, $1<\beta_2, \widetilde{\beta_2}\in\R$, are given constants.\\

$\bullet$ $4<\overline{\beta}_1\in\R$, $1<90\overline{\beta_1}^2/(\overline{\beta_1}^2-16)<\overline{\beta_2}\in\R$, are given constants.\\

$\bullet$ $\alpha_1(t)\equiv\gamma(t)+\beta_1c^2_1(t)+\beta_2c^2_2(t)$,\, $\alpha_2(t)\equiv \overline{\gamma(t)}+\overline{\beta_1} c_1(t)+\overline{\beta_2}c_2^2(t)$, and $\widetilde{\alpha}(t)\equiv \widetilde{\gamma(t)}+\widetilde{\beta_1}c^2_1(t)+\widetilde{\beta_2}c^2_2(t)$ are assumed positive.\\

$\bullet$ $p_1(t)\equiv\exp\left[\int_0^t\alpha_1(s)ds\right]$,\, $p_2(t)\equiv\exp\left[\int_0^t\alpha_2(s)ds\right]$ and $\widetilde{p}(t)\equiv\exp\left[\int_0^t\widetilde{\alpha}(s)ds\right]$.\\

$\bullet$ $L_{\mathscr{F}}(0,T; \R^d)$ is the space of $\mathscr{F}_t$-progressively measurable $\R^d$-valued processes $\varphi(\cdot)$ such that $\E\int^T_0 |\varphi(t)|dt<\infty$.\\

$\bullet$ $M^2(\Omega,\mathscr{F}_T,\P;\R^d)$ is the space of all $\mathscr{F}_T$-measurable $\R^d$-valued random variables $\zeta$ such that
$\E[|\zeta|^2]<\infty$.\\

$\bullet$ $M^2(0,T; \R^d)$ is the space of $\mathscr{F}_t$-progressively measurable $\R^d$-valued processes $\varphi(\cdot)$ such that
$\E\int_0^T|\varphi(t)|^2dt<\infty$.\\

$\bullet$ $\widehat{M_i}^2(\Omega,\mathscr{F}_T,\P;\R^d)$ (resp. ${\widetilde{M}}^2(\Omega,\mathscr{F}_T,\P;\R^d)$) is the space of all $\mathscr{F}_T$-measurable $\R^d$-valued random variables $\xi$ such that $\E[p_i(T)|\xi|^2]<\infty$ (resp. $\E[\widetilde{p}(T)|\xi|^2]<\infty$), $i=1,2$.\\

$\bullet$ $\widehat{M_i}^2(0,T;\R^d)$ (resp. ${\widetilde{M}}^2(0,T;\R^d)$) is the space of $\mathscr{F}_t$-progressively measurable $\R^d$-valued processes $\varphi(\cdot)$ such that $\E\int_0^Tp_i(t)|\varphi(t)|^2dt<\infty$ (reps. $\widetilde{\|\varphi\|}\equiv\E\int_0^T\widetilde{p}(t)|\varphi(t)|^2dt<\infty$), $i=1,2$.\\

$\bullet$ $\widehat{H_i}^2(0,T;\R^d)$ (resp. $\widetilde{H}^2(0,T;\R^d)$) is the space of c\`adl\`ag $\mathscr{F}_t$-adapted $\R^d$-valued processes $\varphi(\cdot)$ such that $\E\left[\sup_{t\in[0,T]}p_i(t)|\varphi(t)|^2\right]<\infty$ (resp.\\
$\E\left[\sup_{t\in[0,T]}\widetilde{p}(t)|\varphi(t)|^2\right]<\infty$), $i=1,2$.\\

\section{Unbounded Lipschitz-type generator}\label{LIP}
In this section, we give sufficient conditions for the existence and uniqueness of a solution pair for (\ref{bsde}). We say that the progressively measurable function $f$ and the random variable $\xi$, or the pair $(f,\xi)$, satisfies {\it conditions A1} (resp. {\it conditions A2}) if:\\
\\
(i)   $\xi\in \widehat{M}_1^2(\Omega,\mathscr{F}_T,\P;\R^d)$\, (resp. $\xi\in \widehat{M}_2^2(\Omega,\mathscr{F}_T,\P;\R^d)$);\\
\\
(ii) $|f(t,y_1,z_1)-f(t,y_2,z_2)|\leq c_1(t)|y_1-y_2|+c_2(t)|z_1-z_2|$, for all $y_1, y_2\in\R^d$, $z_1, z_2\in\R^{d\times k}$, $(t,\omega)$
$a.e.$;\\
\\
(iii) $\left[f(\cdot,0,0)\alpha_1(\cdot)^{-\frac{1}{2}}\right]\in \widehat{M}_1^2(0,T;\R^d)$\, (resp. $\left[f(\cdot,0,0)\alpha_2(\cdot)^{-\frac{1}{2}}\right]\in \widehat{M}_2^2(0,T;\R^d))$.\\

The sufficient conditions  for the solvability of (\ref{bsde}), as given in~\cite{KH}, are similar to our conditions A2. Indeed, if we choose $\overline{\gamma}(t)=0$, $\overline{\beta}_1=\overline{\beta}_2\equiv \beta$, where $\beta$ is {\it large enough}, then conditions A2 are those of~\cite{KH}. Clearly, due to the process $\overline{\gamma}(t)$ our conditions A2 are more general then those of~\cite{KH}. The importance of this process is that assumption (iii) above can be suitable weakened by choosing large values for this process, which is not an option in~\cite{KH}. Moreover, even if we take $\overline{\gamma}(t)=0$, our assumption (i) is weaker than that of~\cite{KH}. Indeed, the parameter $\beta$ of~\cite{KH} should be bigger than $446.05$ (in~\cite{KH} it is only claimed that this coefficient should be {\it large enough}\footnote{For readers' convenience only, we have included an appendix showing that a straightforward calculation gives this numerical lower bound}). This is clearly not the case in conditions A2 where the coefficient $\overline{\beta}_1$ is only required to be greater than $4$.\\

The conditions A1 are new.  In general, these are not comparable with conditions A2. However, in certain special cases we can compare them. For example, if $c_1(t)=0$, $1<\beta_2<\overline{\beta}_2$, $\gamma(t)=\overline{\gamma}(t)$, then $\widehat{M}_1^2(\Omega,\mathscr{F}_T,\P;\R^d)\subset \widehat{M}_2^2(\Omega,\mathscr{F}_T,\P;\R^d)$, and thus the above assumption (i) on the random variable $\xi$ is weaker in the case of conditions A1. Similarly, if $c_2(t)=0$, $\gamma(t)=\overline{\gamma}(t)$, $\overline{\beta}_1=2\beta_1$, then the above assumption (i) on the random variable $\xi$ is weaker in the case of conditions A2.\\

\subsection{Solvability}
In this section we give sufficient conditions for the existence and uniqueness of a solution pair for (\ref{bsde}). Our method of proof is different from~\cite{KH} being based on Picard iterations, and similarly to~\cite{PP}, we begin with a simpler form of \eqref{bsde} and progress towards the
general case. The proofs of the results under conditions A1 and A2 are different and are thus given separately in most cases, but there are also similarities between them.
\begin{lemma}\label{first} Let $\phi(\cdot)\in \widehat{H}_1^2(0,T;\R^d)$, $\psi(\cdot)\in \widehat{M}_1^2(0,T;\R^{d\times k})$ be given, and assume that $\sqrt{\alpha_1(\cdot)}\phi(\cdot)\in \widehat{M}_1^2(0,T;\R^d)$. If the pair $(f,\xi)$ satisfies the conditions A1, then:\\
\\
(i) there exists a unique solution pair $(y(\cdot),z(\cdot))\in \widehat{H}_1^2(0,T;\R^d)\times\widehat{M}_1^2(0,T;\R^{d\times k})$ of equation
\begin{equation}
\label{bsde11}
y(t)= \xi+ \int_t^T f(s, \phi(s), \psi(s))ds-\int_t^T z(s)dW(s),\quad t\in[0,T],
\end{equation}
and $\sqrt{\alpha_1(\cdot)}y(\cdot)\in\widehat{M}_1^2(0,T;\R^d)$.\\
\\
(ii) if $y^+(t)\equiv\mathds{1}_{[y(t)>0]}y(t)$, the processes
\begin{eqnarray}
\int_t^Tp_1(s)y(s)z(s)dW(s) \quad\mbox{{\it and}}\quad\int_t^Tp_1(s)y^+(s)z(s)dW(s),\nonumber
\end{eqnarray}
are martingales.
\end{lemma}
\proof (i) By making use of the Cauchy-Schwartz inequality, we first show that $\int_0^T f(s, \phi(s), \psi(s))ds$ belongs to $M^2(\Omega,\mathscr{F}_T,\P;\R^d)$:
\begin{equation}
\label{estimates}
\begin{split}
&\E\,\bigg|\int_0^T f(s, \phi(s), \psi(s))ds\bigg|^2=\,\E\,\bigg|\int_0^T \sqrt{p_1^{-1}(s)\alpha_1(s)} \frac{\sqrt{p_1(s)}f(s, \phi(s),
\psi(s))}{\sqrt{\alpha_1(s)}}ds\bigg|^2\\
\\
\leq&\,\E\left\{\left[\int_0^Tp_1^{-1}(s)\alpha_1(s)ds\right]\left[\int_0^T\frac{p_1(s)|f(s, \phi(s), \psi(s))|^2}{\alpha_1(s)}ds\right]\right\}\\
\\
\leq&\, \E\int_0^T\frac{p_1(s)|f(s, \phi(s), \psi(s))|^2}{\alpha_1(s)}ds
\end{split}
\end{equation}

\begin{equation}
\begin{split}
=&\,\E\int_0^T\frac{p_1(s)}{\alpha_1(s)}|f(s,\phi(s),\psi(s)-f(s,0,0)+f(s,0,0)|^2ds\nonumber\\
\nonumber\\
\leq&\, \E\int_0^T\frac{p_1(s)}{\alpha_1(s)}[|f(s,\phi(s),\psi(s)-f(s,0,0)|+|f(s,0,0)|]^2ds\nonumber\\
\nonumber\\
\leq&\, \E\int_0^T\frac{p_1(s)}{\alpha_1(s)}[c_1(s)|\phi(s)|+c_2(s)|\psi(s)|+|f(s,0,0)|]^2ds\nonumber\\
\nonumber\\
\leq&\, \E\int_0^T\frac{p_1(s)}{\alpha_1(s)}[3c^2_1(s)|\phi(s)|^2+3c^2_2(s)|\psi(s)|^2+3|f(s,0,0)|^2]ds\nonumber\\
\nonumber\\
=&\, \E\int_0^T\frac{3p_1(s)}{\beta_1}\frac{\beta_1c^2_1(s)}{\gamma(s)+\beta_1c_1^2(s)+\beta_2c_2^2(s)}|\phi(s)|^2ds+\frac{3p_1(s)}{\beta_2}\frac{\beta_2c^2_2(s)}{\gamma(s)+\beta_1c_1^2(s)+\beta_2c_2^2(s)}|\psi(s)|^2ds\nonumber\\
\nonumber\\
&+3\,\E\int_0^T\frac{p_1(s)|f(s,0,0)|^2}{\alpha_1(s)}ds\nonumber\\
\nonumber\\
\leq&\, \frac{3}{\beta_1}\E\int_0^Tp_1(s)|\phi(s)|^2ds+\frac{3}{\beta_2}\E\int_0^Tp_1(s)|\psi(s)|^2ds+3\,\E\int_0^T\frac{p_1(s)|f(s,0,0)|^2}{\alpha_1(s)}ds<\infty\nonumber\\
\nonumber
\end{split}
\end{equation}
Since $\xi\in \widehat{M}_1^2(\Omega,\mathscr{F}_T,\P;\R^d)$ implies that $\xi\in M^2(\Omega,\mathscr{F}_T,\P;\R^d)$, it follows from Lemma 2.1
of~\cite{PP} that \eqref{bsde11} has a unique solution pair $(y(\cdot),z(\cdot))\in M^2(0,T;\R^d)\times M^2(0,T;\R^{d\times k})$. Moreover, since we
proved that \eqref{estimates} is finite, it follows from Lemma 6.2 \footnote{Note that the results in Lemma 6.2 of \cite{KH} is valid for any $\alpha(t)$ (in the notion of \cite{KH}).} of~\cite{KH} that in fact $(y(\cdot),z(\cdot))\in \widehat{H}_1^2(0,T;\R^d)\times\widehat{M}_1^2(0,T;\R^{d\times k})$ and $[\sqrt{\alpha_1(\cdot)}y(\cdot)]\in\widehat{M}_1^2(0,T;\R^d)$.\\

(ii) The proof follows closely that in~\cite{ZJ} (pp. 307), and since it is short, we include it here for completeness. From the Burkholder-Davis-Gundy inequality (see, for example, Theorem 1.5.4 in~\cite{XYZ}), there exists a constant $K$ such that
\begin{eqnarray}
\E\left[\sup_{t\in[0,T]}\left|\int_0^tp_1(s)y(s)z(s)dW(s)\right|\right]\leq K\,\E\left[\int_0^T|\sqrt{p_1(s)}y(s)|^2|\sqrt{p_1(s)}z(s)|^2ds\right]^{\frac{1}{2}}\nonumber\\
\nonumber\\
\leq K\,\E\left[\sup_{t\in [0,T]}|\sqrt{p_1(t)}y(t)|^2\int_0^T\sqrt{p_1(s)}z(s)|^2ds\right]^{\frac{1}{2}}\nonumber\\
\leq \frac{K}{2}\E\left[\sup_{t\in[0,T]}|\sqrt{p_1(t)}y(t)|^2+\int_0^T|\sqrt{p_1(s)}z(s)|^2ds\right]<\infty,\nonumber
\end{eqnarray}
where the last step follows from the fact that $y(\cdot)\in\widehat{H}_1^2(0,T;\R^d)$, $z(\cdot)\in\widehat{M}_1^2(0,T;\R^{d\times k})$, proved in part (i). The conclusion then follows from Corollary 7.22 of~\cite{K}. Since $\sup_{t\in[0,T]}|\sqrt{p_1(t)}y^+(t)|^2\leq\sup_{t\in[0,T]}|\sqrt{p_1(t)}y(t)|^2$,
the conclusion follows even for $\int_0^tp_1(s)y^+(s)z(s)dW(s)$.\qed

\begin{lemma}\label{first2} Let $\phi(\cdot)\in \widehat{H}_2^2(0,T;\R^d)$, $\psi(\cdot)\in \widehat{M}_2^2(0,T;\R^{d\times k})$ be given, and assume that $\sqrt{\alpha_2(\cdot)}\phi(\cdot)\in \widehat{M}_2^2(0,T;\R^d)$. If the pair $(f,\xi)$ satisfies the conditions A2, then:\\
\\
(i) there exists a unique solution pair $(y(\cdot),z(\cdot))\in \widehat{H}_2^2(0,T;\R^d)\times\widehat{M}_2^2(0,T;\R^{d\times k})$ of equation
\begin{equation}
\label{bsde10}
y(t)= \xi+ \int_t^T f(s, \phi(s), \psi(s))ds-\int_t^T z(s)dW(s),\quad t\in[0,T],
\end{equation}
and $\sqrt{\alpha_2(\cdot)}y(\cdot)\in\widehat{M}_2^2(0,T;\R^d)$.\\
\\
(ii) if $y^+(t)\equiv\mathds{1}_{[y(t)>0]}y(t)$, the processes
\begin{eqnarray}
\int_t^Tp_2(s)y(s)z(s)dW(s) \quad\mbox{{\it and}}\quad\int_t^Tp_2(s)y^+(s)z(s)dW(s),\nonumber
\end{eqnarray}
are martingales.
\end{lemma}
\proof The proof of part (ii) is the same as the proof of part (ii) of the previous lemma. We thus focus on part (i). We have
\begin{equation}
\begin{split}
&\E\bigg|\int_0^T f(s, \phi(s), \psi(s))ds\bigg|^2\\
\\
\leq&\, \E\int_0^T\frac{p_2(s)}{\alpha_2(s)}[3c^2_1(s)|\phi(s)|^2+3c^2_2(s)|\psi(s)|^2+3|f(s,0,0)|^2]ds\\
\\
\leq&\,
\E\int_0^T\frac{3p_2(s)}{\overline{\beta_1}^2}\frac{\overline{\beta_1}c_1(s)}{\overline{\gamma(s)}+\overline{\beta_1}c_1(s)+\overline{\beta_2}c_2^2(s)}(\overline{\gamma(s)}+\overline{\beta_1} c_1(s)+\overline{\beta_2}c_2^2(s))|\phi(s)|^2ds\\
\\
&+\,\E\int_0^T\frac{3p_2(s)}{\overline{\beta_2}}\frac{\overline{\beta_2}c^2_2(s)}{\overline{\gamma(s)}+\overline{\beta_1} c_1(s)+\overline{\beta_2}c_2^2(s)}|\psi(s)|^2ds+3\,\E\int_0^T\frac{p_2(s)|f(s,0,0)|^2}{\alpha_2(s)}ds\\
\\
\leq&\, \frac{3}{\overline{\beta_1}^2}\E\int_0^T p_2(s)\alpha_2(s)|\phi(s)|^2ds+\frac{3}{\overline{\beta_2}}\E\int_0^Tp_2(s)|\psi(s)|^2ds+3\,\E\int_0^T\frac{p_2(s)|f(s,0,0)|^2}{\alpha_2(s)}ds\\
\\
<&\,\infty.
\end{split}\nonumber
\end{equation}
The rest of the proof is the same as in the proof of part (i) of the previous lemma.\qed\\

\begin{lemma}\label{second} (i) Let $\phi(\cdot)\in \widehat{H}_1^2(0,T;\R^d)$ be given. If the pair $(f,\xi)$ satisfies conditions A1, then there exists a unique solution pair $(y(\cdot),z(\cdot))\in
\widehat{H}_1^2(0,T;\R^d)\times\widehat{M}_1^2(0,T;\R^{d\times k})$ of equation
\begin{equation}
\label{bsde12}
y(t)= \xi+ \int_t^T f(s, \phi(s), z(s))ds-\int_t^T z(s)dW(s),\quad t\in[0,T],
\end{equation}
and $\sqrt{\alpha_1(\cdot)}y(\cdot)\in\widehat{M}_1^2(0,T;\R^d)$.\\
\\
(ii) Let $\phi(\cdot)\in \widehat{H}_2^2(0,T;\R^d)$ be given. If the pair $(f,\xi)$ satisfies conditions A2, then there exists a unique solution pair $(y(\cdot),z(\cdot))\in
\widehat{H}_2^2(0,T;\R^d)\times\widehat{M}_2^2(0,T;\R^{d\times k})$ of equation
\begin{equation}
y(t)= \xi+ \int_t^T f(s, \phi(s), z(s))ds-\int_t^T z(s)dW(s),\quad t\in[0,T],\nonumber
\end{equation}
and $\sqrt{\alpha_2(\cdot)}y(\cdot)\in\widehat{M}_2^2(0,T;\R^d)$.\\
\end{lemma}
\proof  (i) ({\it Uniqueness})  Let $(y_1(\cdot),z_1(\cdot))$ and $(y_2(\cdot),z_2(\cdot))$ be two solution pairs of \eqref{bsde12} with the claimed properties.
Then
\begin{equation}
\label{pr11}
\begin{split}
&-d\,p_1(t)|y_1(t)-y_2(t)|^2\\
\\
=& \{-\alpha_1(t)p_1(t)\,|y_1(t)-y_2(t)|^2+2p_1(t)(y_1(t)-y_2(t))'\left[f(t,\phi(t),z_1(t))-f(t,\phi(t),z_2(t))\right]\\
\\
&-p_1(t)|z_1(t)-z_2(t)|^2\}dt- 2p_1(t)(y_1(t)-y_2(t))'(z_1(t)-z_2(t))dW(t)
\end{split}
\end{equation}
By using the Lipschitz property of $f$, we have
\begin{equation}
\begin{split}
&-d\,p_1(t)|y_1(t)-y_2(t)|^2\\
\\
\leq& [-\alpha_1(t)p_1(t)\,|y_1(t)-y_2(t)|^2-p_1(t)|z_1(t)-z_2(t)|^2]dt- 2p_1(t)(y_1(t)-y_2(t))'(z_1(t)-z_2(t))dW(t)\\
\\
&+2p_1(t)|y_1(t)-y_2(t)||f(t,\phi(t),z_1(t))-f(t,\phi(t),z_2(t))|dt\\
\\
\leq& [-\alpha_1(t)p_1(t)\,|y_1(t)-y_2(t)|^2-p_1(t)|z_1(t)-z_2(t)|^2]dt- 2p_1(t)(y_1(t)-y_2(t))'(z_1(t)-z_2(t))dW(t)\\
\\
&+2p_1(t)c_2(t)|y_1(t)-y_2(t)||z_1(t)-z_2(t)|dt\\
\\
\leq&[-\alpha_1(t)p_1(t)\,|y_1(t)-y_2(t)|^2-p_1(t)|z_1(t)-z_2(t)|^2]dt- 2p_1(t)(y_1(t)-y_2(t))'(z_1(t)-z_2(t))dW(t)\\
\\
&+\beta_2 c_2^2(t)p_1(t)|y_1(t)-y_2(t)|^2dt+\beta_2^{-1}p_1(t)|z_1(t)-z_2(t)|^2dt\\
\\
\leq&-2p_1(t)(y_1(t)-y_2(t))'(z_1(t)-z_2(t))dW(t),\nonumber
\end{split}
\end{equation}
which in integral form becomes
\begin{eqnarray}
p_1(t)|y_1(t)-y_2(t)|^2\leq \int_t^T-2p_1(s)(y_1(s)-y_2(s))'(z_1(s)-z_2(s))dW(s).\label{super1}
\end{eqnarray}
The stochastic integral in \eqref{super1} is a local martingale that is clearly lower bounded by zero, and is thus a supermartingale (see, for example,
Theorem  7.23 of~\cite{K}). Taking the expectation of both sides of \eqref{super1} results in
\begin{equation}
\begin{split}
\E\left[p_1(t)|y_1(t)-y_2(t)|^2\right]\leq& - \E\left[\int^T_t2p_1(s)(y_1(s)-y_2(s))'(z_1(s)-z_2(s))dW(s)\right]\leq 0.\nonumber
\end{split}
\end{equation}
Since $p_1(t)>0$, it follows that $y_1(t)=y_2(t)$, $\forall$ $t\in[0,T]$, a.s., which proves the uniqueness of $y(\cdot)$. Due to this fact, the integral
form of \eqref{pr11} becomes
\[
0=\int_t^Tp_1(s)|z_1(s)-z_2(s)|^2ds,
\]
which implies that $z_1(t)=z_2(t)$ for $a.e.$ $t\in[0,T]$, and thus proves the uniqueness of $z(\cdot)$.\\

({\it Existence}) Let $z_0(t)\equiv 0$, $\forall t\in[0,T]$,  and for $n\geq 1$ consider the following sequence of equations:
\begin{equation}
\label{bsde1}
y_n(t)= \xi+ \int_t^T f(s, \phi(s), z_{n-1}(s))ds-\int_t^T z_n(s)dW(s),\quad t\in[0,T].
\end{equation}
From Lemma \ref{first} we know that these equations have unique solution pairs $\{(y_n(\cdot),z_n(\cdot))\in \widehat{H}_1^2(0,T; \R^d)\times
\widehat{M}_1^2(0,T; \R^{d\times k})\}_{n\geq 1}$, for which it also holds that $\{\sqrt{\alpha_1(\cdot)}y_n(\cdot)\in\widehat{M}_1^2(0,T;\R^d)\}_{n\geq1}$.
Similarly to the proof of uniqueness, we have
\begin{equation}
\begin{split}
&-d\,p_1(t)|y_{n+1}(t)-y_n(t)|^2\\
\\
=& \{-\alpha_1(t)p_1(t)|y_{n+1}(t)-y_n(t)|^2+2p_1(t)(y_{n+1}(t)-y_n(t))'\left[f(t,\phi(t),z_{n}(t))-f(t,\phi(t),z_{n-1}(t))\right]\\
\\
&-p_1(t)|z_{n+1}(t)-z_n(t)|^2\}dt- 2p_1(t)(y_{n+1}(t)-y_n(t))'(z_{n+1}(t)-z_n(t))dW(t)\\
\\
\leq&[-\alpha_1(t)p_1(t)|y_{n+1}(t)\!-\!y_n(t)|^2\!-\!p_1(t)|z_{n+1}(t)\!-\!z_n(t)|^2]dt\\
\\
&-2p_1(t)(y_{n+1}(t)\!-\!y_n(t))'(z_{n+1}(t)\!-\!z_n(t))dW(t)\\
\\
&+2p_1(t)|y_{n+1}(t)-y_n(t)|\left|f(t,\phi(t),z_{n}(t))-f(t,\phi(t),z_{n-1}(t))\right|dt\\
\\
\leq&[-\alpha_1(t)p_1(t)|y_{n+1}(t)\!-\!y_n(t)|^2\!-\!p_1(t)|z_{n+1}(t)\!-\!z_n(t)|^2]dt\\
\\
&-2p_1(t)(y_{n+1}(t)\!-\!y_n(t))'(z_{n+1}(t)\!-\!z_n(t))dW(t)\\
\\
&+2p_1(t)c_2(t)|y_{n+1}(t)-y_n(t)||z_{n}(t)-z_{n-1}(t)|dt\\
\\
\leq&[-\alpha_1(t)p_1(t)|y_{n+1}(t)\!-\!y_n(t)|^2\!-\!p_1(t)|z_{n+1}(t)\!-\!z_n(t)|^2]dt\\
\\
&-2p_1(t)(y_{n+1}(t)\!-\!y_n(t))'(z_{n+1}(t)\!-\!z_n(t))dW(t)\\
\\
&+\beta_2 c_2^2(t)p_1(t)|y_{n+1}(t)-y_n(t)|^2dt+\beta_2^{-1}p_1(t)|z_{n}(t)-z_{n-1}(t)|^2dt\\
\\
\leq&[-p_1(t)|z_{n+1}(t)\!-\!z_n(t)|^2+\beta_2^{-1}p_1(t)|z_{n}(t)\!-\!z_{n-1}(t)|^2]dt\\
\\
&-\!2p_1(t)(y_{n+1}(t)\!-\!y_n(t))'(z_{n+1}(t)\!-\!z_n(t))dW(t),\nonumber
\end{split}
\end{equation}
which in integral form becomes
\begin{eqnarray}
&&p_1(t)|y_{n+1}(t)-y_n(t)|^2+\int_t^Tp_1(s)|z_{n+1}(s)-z_n(s)|^2ds\nonumber\\
\nonumber\\
&&\leq\beta_2^{-1}\int_t^Tp_1(s)|z_{n}(s)\!-\!z_{n-1}(s)|^2]ds\!-\!\int_t^T2p_1(s)(y_{n+1}(s)\!-\!y_n(s))'(z_{n+1}(s)\!-\!z_n(s))dW(s).\nonumber
\end{eqnarray}
From Lemma \ref{first} (ii), it is clear that the stochastic integral on the right hand side is a martingale. Taking the expected values of both sides gives
\begin{eqnarray}
\E [p_1(t)|y_{n+1}(t)-y_n(t)|^2]+\E\int_t^Tp_1(s)|z_{n+1}(s)-z_n(s)|^2ds\leq\beta_2^{-1}\E\int_t^Tp_1(s)|z_{n}(s)\!-\!z_{n-1}(s)|^2]ds.\nonumber
\end{eqnarray}
Let us define $\eta_n(t)\equiv\E\int_t^Tp_1(s)|y_{n}(s)-y_{n-1}(s)|^2ds$ and $\mu_n(t)\equiv\E\int_t^Tp_1(s)|z_{n}(s)-z_{n-1}(s)|^2ds$. Using the same
argument as in the last part of the proof of Proposition 2.2 in~\cite{PP}, we obtain $\eta_{n+1}(0)\leq\beta_2^{-n}\E\int_0^Tp_1(s)|z_1(s)|^2ds$ and
$\mu_n(0)\leq\beta_2^{-n}\mu_1(0)$. Since the right-hand sides of these two inequalities decrease with $n$, it follows that $\{y_n\}_{n\geq1}$ is a Cauchy sequence in $\widehat{M}^2_1(0,T;\R^d)$, and $\{z_n\}_{n\geq 1}$ is a Cauchy sequence in $\widehat{M}^2_1(0,T;\R^{d\times k})$. Moreover, this also implies that $\{\sqrt{\alpha_1}y_n\}_{n\geq 1}$ is a Cauchy sequence in $\widehat{M}^2_1(0,T;\R^d)$. Hence, the limiting processes $y^*=\lim_{n\rightarrow\infty}y_n$ and $z^*=\lim_{n\rightarrow\infty}z_n$ are the solution pair of \eqref{bsde12}. In addition, when such a pair of processes is substituted in \eqref{bsde12}, then \eqref{bsde12} becomes an example of \eqref{bsde11} with $\psi(\cdot)=z^*(\cdot)$. Therefore, Lemma \ref{first} applies, and we have that $y^*(\cdot)\in \widehat{H}^2_1(0,T;\R^{d\times k})$.\\
\\
(ii) Due to Lemma \ref{first2}, the proof in this case is identical to the proof of part (i) (with an obvious change of notation), and is thus omitted.\qed\\

\begin{theorem}
\label{theorem1}
(i) If the pair $(f,\xi)$ satisfies conditions A1, then equation \eqref{bsde} has a unique solution pair $(y(\cdot), z(\cdot))\in \widehat{H}_1^2(0,T; \R^d)\times \widehat{M}_1^2(0,T; \R^{d\times k})$, and $\sqrt{\alpha_1(\cdot)}y(\cdot)\in\widehat{M}_1^2(0,T; \R^d)$.\\
\\
(ii) If the pair $(f,\xi)$ satisfies conditions A2, then equation \eqref{bsde} has a unique solution pair $(y(\cdot), z(\cdot))\in \widehat{H}_2^2(0,T; \R^d)\times \widehat{M}_2^2(0,T; \R^{d\times k})$, and $\sqrt{\alpha_2(\cdot)}y(\cdot)\in\widehat{M}_2^2(0,T; \R^d)$.
\end{theorem}
\proof
(i) ({\it Uniqueness}) Let $(y_1(\cdot),z_1(\cdot))$ and $(y_2(\cdot),z_2(\cdot))$ be two solution pairs of \eqref{bsde} with the claimed properties. Then
we have
\begin{equation}
\begin{split}
&-d\,p_1(t)|y_1(t)-y_2(t)|^2\\
\\
=& \{-\alpha_1(t)p_1(t)\,|y_1(t)-y_2(t)|^2+2p_1(t)(y_1(t)-y_2(t))'\left[f(t,y_1(t),z_1(t))-f(t,y_2(t),z_2(t))\right]\\
\\
&-p_1(t)|z_1(t)-z_2(t)|^2\}dt- 2p_1(t)(y_1(t)-y_2(t))'(z_1(t)-z_2(t))dW(t)\\
\\
\leq&[-\alpha_1(t)p_1(t)\,|y_1(t)-y_2(t)|^2-p_1(t)|z_1(t)-z_2(t)|^2]dt-2p_1(t)(y_1(t)-y_2(t))'(z_1(t)-z_2(t))dW(t)\\
\\
&+2p_1(t)|y_1(t)-y_2(t)|\left|f(t,y_1(t),z_1(t))-f(t,y_2(t),z_2(t))\right|dt\\
\\
\leq&[-\alpha_1(t)p_1(t)\,|y_1(t)-y_2(t)|^2-p_1(t)|z_1(t)-z_2(t)|^2]dt-2p_1(t)(y_1(t)-y_2(t))'(z_1(t)-z_2(t))dW(t)\\
\\
&+2p_1(t)|y_1(t)-y_2(t)|\left[c_1(t)|y_1(t)-y_2(t)|+c_2(t)|z_1(t)-z_2(t)|\right]dt\\
\\
\leq&[-\alpha_1(t)p_1(t)\,|y_1(t)-y_2(t)|^2-p_1(t)|z_1(t)-z_2(t)|^2]dt-2p_1(t)(y_1(t)-y_2(t))'(z_1(t)-z_2(t))dW(t)\\
\\
&+\beta_1c_1^2p_1(t)|y_1(t)-y_2(t)|^2dt+\beta_1^{-1}p_1(t)|y_1(t)-y_2(t)|^2dt\\
\\
&+\beta_2c_2^2p_1(t)|y_1(t)-y_2(t)|^2dt+\beta_2^{-1}p_1(t)|z_1(t)-z_2(t)|^2dt\\
\\
\leq&[\beta_1^{-1}p_1(t)|y_1(t)-y_2(t)|^2+(\beta_2^{-1}-1)p_1(t)|z_1(t)-z_2(t)|^2]dt\\
\\
&-2p_1(t)(y_1(t)-y_2(t))'(z_1(t)-z_2(t))dW(t)\\
\\
\leq&\beta_1^{-1}p_1(t)|y_1(t)-y_2(t)|^2dt-2p_1(t)(y_1(t)-y_2(t))'(z_1(t)-z_2(t))dW(t).\nonumber\\
\end{split}
\end{equation}
With the help of Lemma \ref{first} (ii) and the Gronwall's lemma, the conclusion follows similarly to the proof of uniqueness in Lemma \ref{second}.\\

({\it Existence})  Let $y_0(t)\equiv 0$, $\forall t\in[0,T]$, and for $n\geq 1$ consider the sequence of equations:
\begin{equation}
\label{bsde1}
y_n(t)= \xi+ \int_t^T f(s, y_{n-1}(s), z_n(s))ds-\int_t^T z_n(s)dW(s),\quad t\in[0,T].
\end{equation}
From Lemma \ref{second} we know that these equations have unique solution pairs $\{(y_n(\cdot),z_n(\cdot))\in \widehat{H}_1^2(0,T; \R^d)\times
\widehat{M}_1^2(0,T; \R^{d\times k})\}_{n\geq 1}$. Then
\begin{equation}
\label{pr}
\begin{split}
&-d\,p_1(t)|y_{n+1}(t)-y_n(t)|^2\\
\\
=& \{-\alpha_1(t)p_1(t)|y_{n+1}(t)-y_n(t)|^2+2p_1(t)(y_{n+1}(t)-y_n(t))'\left[f(t,y_n(t),z_{n+1}(t))-f(t,y_{n-1}(t),z_n(t))\right]\\
\\
&-p_1(t)|z_{n+1}(t)-z_n(t)|^2\}dt- 2p_1(t)(y_{n+1}(t)-y_n(t))'(z_{n+1}(t)-z_n(t))dW(t).\\
\\
\leq&[-\alpha_1(t)p_1(t)|y_{n+1}(t)\!-\!y_n(t)|^2\!-\!p_1(t)|z_{n+1}(t)\!-\!z_n(t)|^2]dt\\
\\
&-2p_1(t)(y_{n+1}(t)\!-\!y_n(t))'(z_{n+1}(t)\!-\!z_n(t))dW(t)\\
\\
&+2p_1(t)|y_{n+1}(t)-y_n(t)|\left|f(t,y_{n}(t),z_{n+1}(t))-f(t,y_{n-1}(t),z_n(t))\right|dt\\
\\
\leq&[-\alpha_1(t)p_1(t)|y_{n+1}(t)\!-\!y_n(t)|^2\!-\!p_1(t)|z_{n+1}(t)\!-\!z_n(t)|^2]dt\\
\\
&-2p_1(t)(y_{n+1}(t)\!-\!y_n(t))'(z_{n+1}(t)\!-\!z_n(t))dW(t)\\
\\
&+2p_1(t)|y_{n+1}(t)-y_n(t)|\left[c_1(t)|y_n(t)-y_{n-1}(t)|+c_2(t)|z_{n+1}(t)-z_n(t)|\right]dt\nonumber\\
\\
\leq&[-\alpha_1(t)p_1(t)|y_{n+1}(t)\!-\!y_n(t)|^2\!-\!p_1(t)|z_{n+1}(t)\!-\!z_n(t)|^2]dt\\
\\
&-2p_1(t)(y_{n+1}(t)\!-\!y_n(t))'(z_{n+1}(t)\!-\!z_n(t))dW(t)\\
\\
&+\beta_1c_1^2(t)p_1(t)|y_{n+1}(t)-y_n(t)|^2dt+\beta_1^{-1}p_1(t)|y_n(t)-y_{n-1}(t)|^2dt\\
\\
&+\beta_2c_2^2(t)p_1(t)|y_{n+1}(t)-y_n(t)|^2dt+ \beta_2^{-1}p_1(t)|z_{n+1}(t)-z_n(t)|^2dt\nonumber\\
\end{split}
\end{equation}
\begin{equation}
\begin{split}
\leq&\,\beta_1^{-1}p_1(t)|y_n(t)-y_{n-1}(t)|^2dt+(\beta_2^{-1}-1)p_1(t)|z_{n+1}(t)-z_n(t)|^2dt\\
\\
&-2p_1(t)(y_{n+1}(t)-y_n(t))'(z_{n+1}(t)-z_n(t))dW(t).\nonumber
\end{split}
\end{equation}
Due to Lemma \ref{first} (ii), the expectation of the integral-form of this inequity becomes
\begin{equation}
\begin{split}
\E\left[p_1(t)|y_{n+1}(t)-y_n(t)|^2\right]\leq&\, \E\int^T_t \beta_1^{-1}p_1(s)|y_n(s)-y_{n-1}(s)|^2ds\\
&+\E\int^T_t\!\!\!(\beta_2^{-1}\!-\!1)p_1(s)|z_{n+1}(s)\!-\!z_n(s)|^2ds,\label{N}
\end{split}
\end{equation}
Using the notation $\nu_{n+1}(t)\equiv\E\int_t^Tp_1(t)|y_{n+1}(s)-y_n(s)|^2ds$, and similarly to the last part of the proof of Theorem 3.1 of~\cite{PP}, we obtain $\nu_{n+1}(0)\leq \beta_1^{-n}\frac{1}{n!}\nu_1(0)$. Since the sum of the right-hand side of this inequality converges, we conclude, together with \eqref{N}, that $\{y_n\}$ is a Cauchy sequence in $\widehat{M}_1^2(0,T; \R^d)$, and $\{z_n\}$ is a Cauchy sequence in $\widehat{M}_1^2(0,T; \R^{d\times k})$. Moreover, this also implies that $\{\sqrt{\alpha}y_n\}_{n\geq1}$ is a Cauchy sequence in $\widehat{M}_1^2(0,T; \R^d)$.
Thus the limiting processes $y^*=\lim_{n\rightarrow\infty}y_n$ and $z^*=\lim_{n\rightarrow\infty}z_n$ are the solution pair to \eqref{bsde}. In
addition, when such a pair of processes is substituted in \eqref{bsde}, then \eqref{bsde} becomes an example of \eqref{bsde11} with
$\phi(\cdot)=y^*(\cdot)$ and $\psi(\cdot)=z^*(\cdot)$. Therefore, Lemma \ref{first} applies, and we have that $y^*(\cdot)\in
\widehat{H}^2_1(0,T;\R^{d\times k})$.\\

(ii) {\it (Uniqueness)} Let $(y_1(\cdot),z_1(\cdot))$ and $(y_2(\cdot),z_2(\cdot))$ be two solution pairs of \eqref{bsde} with the claimed properties. Similarly to the proof of uniqueness for part (i), we have
\begin{equation}
\begin{split}
&-d\,p_2(t)|y_1(t)-y_2(t)|^2\\
\\
\leq&[-\alpha_2(t)p_2(t)\,|y_1(t)-y_2(t)|^2-p_2(t)|z_1(t)-z_2(t)|^2]dt-2p_2(t)(y_1(t)-y_2(t))'(z_1(t)-z_2(t))dW(t)\\
\\
&+2p_2(t)|y_1(t)-y_2(t)|\left[c_1(t)|y_1(t)-y_2(t)|+c_2(t)|z_1(t)-z_2(t)|\right]dt\\
\\
\leq&[-\alpha_2(t)p_2(t)\,|y_1(t)-y_2(t)|^2-p_2(t)|z_1(t)-z_2(t)|^2]dt-2p_2(t)(y_1(t)-y_2(t))'(z_1(t)-z_2(t))dW(t)\\
\\
&+2c_1p_2(t)|y_1(t)-y_2(t)|^2dt+\overline{\beta_2}c_2^2p_2(t)|y_1(t)-y_2(t)|^2dt+\overline{\beta_2}^{-1}p_2(t)|z_1(t)-z_2(t)|^2dt\\
\\
\leq&(\beta_2^{-1}-1)p_2(t)|z_1(t)-z_2(t)|^2dt-2p_2(t)(y_1(t)-y_2(t))'(z_1(t)-z_2(t))dW(t)\\
\\
\leq& -2p_2(t)(y_1(t)-y_2(t))'(z_1(t)-z_2(t))dW(t).\nonumber\\
\end{split}
\end{equation}
Then the expectation of integral-form of this inequality becomes
\[
\E\left[p_2(t)|y_1(t)-y_2(t)|^2\right]\leq \E\left[\int^T_t -2p_2(s)(y_1(s)-y_2(s))'(z_1(s)-z_2(s))dW(s)\right].
\]
Since the right-hand side is a martingale by Lemma \ref{first} (ii), the conclusion follows.\\

{\it (Existence)} Let $y_0(t)\equiv 0$, $\forall t\in[0,T]$, and for $n\geq 1$ consider the sequence of equations:
\begin{equation}
\label{bsde1}
y_n(t)= \xi+ \int_t^T f(s, y_{n-1}(s), z_n(s))ds-\int_t^T z_n(s)dW(s),\quad t\in[0,T].
\end{equation}
From Lemma \ref{second} we know that these equations have unique solution pairs $\{(y_n(\cdot),z_n(\cdot))\in \widehat{H}_2^2(0,T; \R^d)\times
\widehat{M}_2^2(0,T; \R^{d\times k})\}_{n\geq 1}$. By Lemma 6.2 of \cite{KH}, we have following estimates:
\begin{equation}
\label{est1}
\begin{split}
&\E\left[\int^T_0p_2(t)|y_{n+1}(t)-y_n(t)|^2\alpha_2(t)dt\right]\\
\leq&\, 8\,\E\left[\int^T_0p_2(t)\frac{|f(t,y_n(t),z_{n+1}(t))-f(t,y_{n-1}(t),z_{n}(t))|^2}{\alpha_2(t)}dt\right],\\
\end{split}
\end{equation}
and
\begin{equation}
\label{est2}
\begin{split}
&\E\left[\int^T_0p_2(t)|z_{n+1}(t)-z_n(t)|^2dt\right]\\
\leq&\, 45\,\E\left[\int^T_0p_2(t)\frac{|f(t,y_n(t),z_{n+1}(t))-f(t,y_{n-1}(t),z_{n}(t))|^2}{\alpha_2(t)}dt\right].\\
\end{split}
\end{equation}
\\
By the Lipschitz condition, we have
\begin{equation}
\begin{split}
&\E\left[\int^T_0p_2(t)\frac{|f(t,y_n(t),z_{n+1}(t))-f(t,y_{n-1}(t),z_{n}(t))|^2}{\alpha_2(t)}dt\right]\\
\\
\leq&\, \E\int_0^T\frac{p_2(t)}{\alpha_2(t)}[c_1(t)|y_n(t)-y_{n-1}(t)|+c_2(t)|z_{n+1}(t)-z_{n}(t)|]^2dt\\
\\
\leq&\, 2\E\,\int_0^T\frac{p_2(t)}{\alpha_2(t)}[c^2_1(t)|y_n(t)-y_{n-1}(t)|^2+c^2_2(t)|z_{n+1}(t)-z_{n}(t)|^2]dt\\
\\
\leq&\, 2\E\,\int_0^T\frac{p_2(t)}{\overline{\beta_1}^2}\frac{\overline{\beta_1}c_1(t)}{\overline{\gamma(t)}+\overline{\beta_1}c_1(t)+\overline{\beta_2}c_2^2(t)}(\overline{\gamma(t)}+\overline{\beta_1} c_1(t)+\overline{\beta_2}c_2^2(t))|y_n(t)-y_{n-1}(t)|^2dt\\
\\
&+2\E\,\int_0^T\frac{p_2(t)}{\overline{\beta_2}}\frac{\overline{\beta_2}c^2_2(t)}{\overline{\gamma(t)}+\overline{\beta_1} c_1(t)+\overline{\beta_2}c_2^2(t)}|z_{n+1}(t)-z_{n}(t)|]^2ds\\
\\
\leq& \frac{2}{\overline{\beta_1}^2}\E\int_0^T p_2(t)\alpha_2(t)|y_n(t)-y_{n-1}(t)|^2dt+\frac{2}{\overline{\beta_2}}\E\int_0^Tp_2(t)|z_{n+1}(t)-z_{n}(t)|]^2dt.\nonumber\\
\end{split}
\end{equation}
\\
Substituting it into \eqref{est2}, we have
\begin{equation}
\begin{split}
&\E\left[\int^T_0p_2(t)|z_{n+1}(t)-z_n(t)|^2dt\right]\\
\\
\leq& \frac{90}{\overline{\beta_1}^2}\E\int_0^T p_2(t)\alpha_2(t)|y_n(t)-y_{n-1}(t)|^2dt+\frac{90}{\overline{\beta_2}}\E\int_0^Tp_2(t)|z_{n+1}(t)-z_{n}(t)|]^2dt.\nonumber
\end{split}
\end{equation}
\\
Let $\overline{\beta_2}>90$. Then we have
\begin{equation}\label{z}
\E\left[\int^T_0p_2(t)|z_{n+1}(t)-z_n(t)|^2dt\right] \leq \frac{\frac{90}{\overline{\beta_1}^2}}{\left(1-\frac{90}{\overline{\beta_2}}\right)}\E\int_0^T p_2(t)\alpha_2(t)|y_n(t)-y_{n-1}(t)|^2dt.
\end{equation}
\\
Substituting it into \eqref{est2}, we obtain
\begin{equation}
\begin{split}
&\E\left[\int^T_0p_2(t)|y_{n+1}(t)-y_n(t)|^2\alpha_2(t)dt\right]\\
\\
\leq& \frac{16}{\overline{\beta_1}^2}\E\int_0^T p_2(t)\alpha_2(t)|y_n(t)-y_{n-1}(t)|^2dt+\frac{16}{\overline{\beta_2}}\E\int_0^Tp_2(t)|z_{n+1}(t)-z_{n}(t)|]^2dt\\
\\
\leq& \frac{16}{\overline{\beta_1}^2}\E\int_0^T p_2(t)\alpha_2(t)|y_n(t)-y_{n-1}(t)|^2dt+ \frac{16}{\overline{\beta_2}} \frac{\frac{90}{\overline{\beta_1}^2}}{\left(1-\frac{90}{\overline{\beta_2}}\right)}\E\int_0^T p_2(t)\alpha_2(t)|y_n(t)-y_{n-1}(t)|^2dt\\
\\
=& \left[\frac{16}{\overline{\beta_1}^2}+\frac{16}{\overline{\beta_2}} \frac{\frac{90}{\overline{\beta_1}^2}}{\left(1-\frac{90}{\overline{\beta_2}}\right)}\right]\E\int_0^T p_2(t)\alpha_2(t)|y_n(t)-y_{n-1}(t)|^2dt.\nonumber
\end{split}
\end{equation}

Let $\kappa=\left[\frac{16}{\overline{\beta_1}^2}+\frac{16}{\overline{\beta_2}} \frac{\frac{90}{\overline{\beta_1}^2}}{\left(1-\frac{90}{\overline{\beta_2}}\right)}\right]$, i.e. $\overline{\beta}_1>4$ and $\overline{\beta}_2>\frac{90\overline{\beta}^2_1}{\overline{\beta}_1^2-16}$. Then the conclusion follows similarly to the proof of existence in Lemma \ref{second}.

\subsection{Comparison theorem}\label{C}
The following results generalise Peng's comparison theorem (\cite{Peng2},~\cite{Peng}) to equations with a possibly unbounded generator. Similarly to~\cite{Peng2},~\cite{Peng}, we assume that $d=1$. In addition to equation (\ref{bsde}), let us consider two further equations
\begin{eqnarray}
{\widehat{y}_1}(t)= {\widehat{\xi}_1}+ \int_t^T [\widehat{f}_1(s,{\widehat{y}_1}(s),{\widehat{z}_1}(s))]ds-\int_t^T {\widehat{z}_1}(s)dW(s),\quad t\in[0,T],\nonumber\\
\nonumber\\
{\widehat{y}_2}(t)= {\widehat{\xi}_2}+ \int_t^T [\widehat{f}_2(s,{\widehat{y}_2}(s),{\widehat{z}_2}(s))]ds-\int_t^T {\widehat{z}_2}(s)dW(s),\quad t\in[0,T].\nonumber
\end{eqnarray}
We assume that the pair $(\widehat{f}_1,\widehat{\xi}_1)$ satisfies conditions A1, whereas the pair $(\widehat{f}_2,\widehat{\xi}_2)$ satisfies conditions A2. Based on Theorem \ref{theorem1}, this means that there exist unique solution pairs $(\widehat{y}_1(\cdot), \widehat{z}_1(\cdot))\in \widehat{H}_1^2(0,T; \R)\times \widehat{M}_1^2(0,T; \R^{1\times k})$ and $(\widehat{y}_2(\cdot), \widehat{z}_2(\cdot))\in \widehat{H}_2^2(0,T; \R)\times \widehat{M}_2^2(0,T; \R^{1\times k})$. The following differences will appear in the proof:
\begin{eqnarray}
Y_1(t)\equiv y(t)-\widehat{y}_1(t),\quad Z_1(t)\equiv z(t)-\widehat{z}_1(t),\nonumber\\
Y_2(t)\equiv y(t)-\widehat{y}_2(t),\quad Z_2(t)\equiv z(t)-\widehat{z}_2(t).\nonumber
\end{eqnarray}
\begin{theorem}
\label{compare} (Comparison theorem) (i) If $\widehat{\xi}_1\geq \xi$ and $\widehat{f}_1(t,y,z)\geq f(t,y,z)$, $a.s.$ $\forall$  $(t,y,z)\in [0,T]\times \R\times \R^{1\times k}$, then ${\widehat{y}_1}(t)\geq y(t)$, $\forall$ $t\in[0,T]$, $a.s.$.\\
\\
(ii) If $\widehat{\xi}_2\geq \xi$ and $\widehat{f}_2(t,y,z)\geq f(t,y,z)$, $a.s.$ $\forall$  $(t,y,z)\in [0,T]\times \R\times \R^{1\times k}$, then ${\widehat{y}_2}(t)\geq y(t)$, $\forall$ $t\in[0,T]$, $a.s.$.
\end{theorem}
\proof (i) The equation of the difference $Y_1(t)$ is
\begin{eqnarray}
-dY_1(t)=[f(t,y(t),z(t))-\widehat{f}_1(t,{\widehat{y}_1}(t),{\widehat{z}_1}(t))]dt-Z_1(t)dW(t).\nonumber
\end{eqnarray}
Denoting by $Y_1^+(t)\equiv \mathds{1}_{[Y_1(t)>0]}Y_1(t)$, and using Tanaka-Meyer formula (see Theorem 6.1.2 in~\cite{RY}), we obtain
\begin{eqnarray}
-dY_1^+(t)=-\mathds{1}_{[Y_1(t)>0]}dY_1(t)-\frac{1}{2}dL(t),\nonumber\
\end{eqnarray}
where $L(t)$ is the local time of $Y_1(\cdot)$ at $0$. Since $\int_0^T|Y_1(t)|dL(t)=0$, a.s. (see Proposition 6.1.3 in~\cite{RY}), we have
\begin{eqnarray}
-d[Y_1^+(t)]^2&=&2Y_1^+(t)\mathds{1}_{[Y_1(t)>0]}[f(t,y(t),z(t))-\widehat{f}_1(t,{\widehat{y}_1}(t),{\widehat{z}_1}(t))]dt\nonumber\\
\nonumber\\
&-&\mathds{1}_{[Y_1(t)>0]}Z_1^2(t)dt-\mathds{1}_{[Y_1(t)>0]}2Y_1^+(t)Z_1(t)dW(t).\nonumber
\end{eqnarray}
Using It\^{o} formula, we obtain
\begin{equation}
\begin{split}
&-d\,p_1(t)[Y_1^+(t)]^2\\
\\
=&-\alpha_1(t)p_1(t)[Y_1^+(t)]^2dt+2p_1(t)Y_1^+(t)\mathds{1}_{[Y_1(t)>0]}[f(t,y(t),z(t))-\widehat{f}_1(t,{\widehat{y}_1}(t),{\widehat{z}_1}(t))]dt\\
\\
&-\mathds{1}_{[Y_1(t)>0]}p_1(t)Z_1^2(t)dt-2p_1(t)Y_1^+(t)Z_1(t)dW(t)\\
\\
\leq& -\!\alpha_1(t)p_1(t)[Y_1^+(t)]^2dt+2p_1(t)Y_1^+(t)\mathds{1}_{[Y_1(t)>0]}[f(t,y(t),z(t))-\widehat{f}_1(t,y(t),z(t))\\
\\
&+\widehat{f}_1(t,y(t),z(t))-\widehat{f}_1(t,{\widehat{y}_1}(t),{\widehat{z}_1}(t))]dt\\
\\
&-\mathds{1}_{[Y_1(t)>0]}p_1(t)Z_1^2(t)dt-2p_1(t)Y_1^+(t)Z_1(t)dW(t)\\
\\
\leq&[-\!\alpha_1(t)p_1(t)[Y_1^+(t)]^2\!+\!2p_1(t)Y_1^+(t)\mathds{1}_{[Y_1(t)>0]}[f(t,y(t),z(t))-\widehat{f}_1(t,y(t),z(t))]\!\\
\\
&-\!\mathds{1}_{[Y_1(t)>0]}p_1(t)Z_1^2(t)]dt\!-\!\mathds{1}_{[Y_1(t)>0]}2p_1(t)Y_1^+(t)Z_1(t)dW(t)+\beta_1p_1(t)c_1^2(t)[Y_1^+(t)]^2dt\\
\\
&+\beta_1^{-1}p_1(t)[Y_1^+(t)]^2dt+\beta_2c_2^2(t)p_1(t)[Y_1^+(t)]^2dt+\beta_2^{-1}\mathds{1}_{[Y_1(t)>0]}p_1(t)Z_1^2(t)dt\\
\\
\leq&\beta_1^{-1}p_1(t)[Y_1^+(t)]^2dt-2p_1(t)Y_1^+(t)Z_1(t)dW(t),\nonumber
\end{split}
\end{equation}
which in integral form becomes
\begin{eqnarray}
p_1(t)[Y_1^+(t)]^2\leq \int_t^T\beta_1^{-1}p_1(s)[Y_1^+(s)]^2ds-\int_t^T 2p_1(s)Y_1^+(s)Z_1(s)dW(s).\nonumber
\end{eqnarray}
The stochastic integral on the right-hand side is a martingale due to Lemma \ref{first} (ii). Therefore,
\begin{eqnarray}
\E[p_1(t)[Y_1^+(t)]^2]\leq \E\int_t^T\beta_1^{-1}p_1(s)[Y_1^+(s)]^2ds,\nonumber
\end{eqnarray}
and the conclusion follows from Gronwall's lemma.\\

(ii) In a similar way to the proof of part (i), we have

\begin{equation}
\begin{split}
&-d\,p_2(t)[Y_2^+(t)]^2\\
\\
=&-\alpha_2(t)p_2(t)[Y_2^+(t)]^2dt+2p_2(t)Y_2^+(t)\mathds{1}_{[Y_2(t)>0]}[f(t,y(t),z(t))-\widehat{f}_2(t,{\widehat{y}_2}(t),{\widehat{z}_2}(t))]dt\\
\\
&-\mathds{1}_{[Y_2(t)>0]}p_2(t)Z_2^2(t)dt-2p_2(t)Y_2^+(t)Z_2(t)dW(t)\\
\\
\leq& -\!\alpha_2(t)p_2(t)[Y_2^+(t)]^2dt+2p_2(t)Y_2^+(t)\mathds{1}_{[Y_2(t)>0]}[f(t,y(t),z(t))-\widehat{f}_2(t,y(t),z(t))\\
\\
&+\widehat{f}_2(t,y(t),z(t))-\widehat{f}_2(t,{\widehat{y}_2}(t),{\widehat{z}_2}(t))]dt\\
\\
&-\mathds{1}_{[Y_2(t)>0]}p_2(t)Z_2^2(t)dt-2p_2(t)Y_2^+(t)Z_2(t)dW(t)\\
\\
\leq&[-\!\alpha_2(t)p_2(t)[Y_2^+(t)]^2\!+\!2p_2(t)Y_2^+(t)\mathds{1}_{[Y_2(t)>0]}[f(t,y(t),z(t))-\widehat{f}_2(t,y(t),z(t))]\!\\
\\
&-\!\mathds{1}_{[Y_2(t)>0]}p_2(t)Z_2^2(t)]dt\!-\!\mathds{1}_{[Y_2(t)>0]}2p_2(t)Y_2^+(t)Z_2(t)dW(t)+2p_2(t)c_1(t)[Y_2^+(t)]^2dt\\
\\
&+\overline{\beta}_2c_2^2(t)p_2(t)[Y_2^+(t)]^2dt+\overline{\beta}_2^{-1}\mathds{1}_{[Y_2(t)>0]}p_2(t)Z_2^2(t)dt\\
\\
\leq&-2\mathds{1}_{[Y_2(t)>0]}p_1(t)Y_2^+(t)Z_2(t)dW(t),\nonumber
\end{split}
\end{equation}
which in integral form becomes
\begin{eqnarray}
p_2(t)[Y_2^+(t)]^2\leq-\int_t^T 2\mathds{1}_{[Y_2(s)>0]}p_2(s)Y_2^+(s)Z_2(s)dW(s).\nonumber
\end{eqnarray}
Since the stochastic integral on the right-hand side is a martingale due to Lemma \ref{first} (ii), we have
\begin{eqnarray}
\E[p_2(t)[Y_2^+(t)]^2]\leq 0,\nonumber
\end{eqnarray}
which concludes the proof.\qed\\

\section{Unbounded continuous generator}\label{CON}

In this section, we consider the one-dimensional version of equation \eqref{bsde} with a continuous $f$ with respect to $y$ and $z$, which satisfies a linear growth condition rather than the Lipschitz-type condition. The main idea here, as in~\cite{LS}, is to approximate the generator $f$ by an infinite sequence of Lipschitz-type approximating functions. Each such a function generates a BSDE, and we show that the solutions to such a sequence of BSDEs converge to the solution of \eqref{bsde}.

We say that the progressively measurable function $f$ and the random variable $\xi$, or the pair $(f,\xi)$, satisfies {\it conditions A3} if:\\
\\
(i)' $d=1$ and $f(t,y,z)$ is a continuous function of $y$ and $z$;\\
\\
(ii)' $|f(t,y,z)|\leq c_0(t)+c_1(t)|y|+c_2(t)|z|$, for all $y\in\R$, $z\in\R^{k}$, $(t,\omega)$  $a.e.$;\\
\\
(iii)'  $\xi\in {{\widetilde{M}}}^2(\Omega,\mathscr{F}_T,\P;\R)$;\\
\\
(iv)' $c_0(\cdot)\in{\widetilde{M}}^2(0,T;\R)$ and $\left[c_0(\cdot)\widetilde{\alpha}(\cdot)^{-\frac{1}{2}}\right]\in{\widetilde{M}}^2(0,T;\R)$.\\

As already mentioned in the introduction, our assumption (ii)' permits for random and possibly unbounded coefficients $c_1(\cdot)$ and $c_2(\cdot)$, which is not the case in \cite{WW} and \cite{WH}. The results of the previous section are our main tools in dealing with the problem of solvability. We focus on utilizing Theorem \ref{theorem1} (i) and Theorem \ref{compare} (i) only, as similar results can be obtained by applying Theorem \ref{theorem1} (ii) and Theorem \ref{compare} (ii).\\

We introduce the sequence of functions
\begin{equation}
f_n(t,y,z)\equiv\sup_{(u,v)\in\R^{1+k}}\{f(t,u,v)-[c_1(t)+n]|u-y|-[c_2(t)+n]|v-z|\},\quad n\geq1,\nonumber
\end{equation}
which are clearly well-defined. Their main properties are summarized in the following result.
\begin{lemma}
\label{seq}
(i) Linear growth: for any $y\in\R$, $z\in \R^k$, $|f_n(t,y,z)|\leq c_0(t)+ c_1(t)|y|+ c_2(t)|z|$;\\
\\
(ii) Monotonicity: $f_n$ is a decreasing function of $n$;\\
\\
(iii) Lipschitz condition: for any $y_1, y_2\in\R$, $z_1, z_2\in \R^k$,
\[
|f_n(t,y_1,z_1)-f_n(t,y_2,z_2)|\leq [c_1(t)+n]|y_1-y_2|+[c_2(t)+n]|z_1-z_2|;
\]
(iv) Convergence: for any $y\in\R$, $z\in \R^k$, $\lim_{n\rightarrow \infty}f_n(t,y,z)=f(t,y,z)$.
\end{lemma}

\proof
(i) By the linear growth of $f$, for all $y\in\R$, $z\in\R^k$, we have
\begin{equation}
\begin{split}
&f_n(t,y,z)\leq \sup_{(u,v)\in\R^{1+k}}\{|f(t,u,v)|-[c_1(t)+n]|u-y|-[c_2(t)+n]|v-z|\}\\
\\
\leq& c_0(t)+\sup_{(u,v)\in\R^{1+k}}\{c_1(t)|u|+ c_2(t)|v|-[c_1(t)+n]|u-y|-[c_2(t)+n]|v-z|\}\\
\\
\leq& c_0(t)+\sup_{(u,v)\in\R^{1+k}}\{c_1(t)|u|+ c_2(t)|v|-c_1(t)(|u|-|y|)-c_2(t)(|v|-|z|)\}\\
\\
\leq& c_1(t)+ c_2(t)|y|+ c_3(t)|z|.\nonumber\\
\end{split}
\end{equation}
The inequality $f_n(t,y,z)\geq- c_0(t)-c_1(t)|y|-c_2(t)|z|$ can be proved similarly.\\

(ii) This follows from the definition of $f_n$ itself.\\

(iii) By inequality $|\sup_{i\in I}a_i-\sup_{i\in I}b_i|\leq \sup_{i\in I}|a_i-b_i|$, with $I$ being an arbitrary index set, we have
\begin{equation}
\begin{split}
&|f_n(t,y_1,z_1)-f_n(t,y_2,z_2)|\\
\\
=& \bigg|\sup_{(u,v)\in\R^{1+k}}\{f(t,u,v)-[c_1(t)+n]|u-y_1|-[c_2(t)+n]|v-z_1|\}\\
\\
& -\sup_{(u,v)\in\R^{1+k}}\{f(t,u,v)-[c_1(t)+n]|u-y_2|-[c_2(t)+n]|v-z_2|\}\bigg|\\
\\
\leq& \sup_{(u,v)\in\R^{1+k}}\big|\{[c_1(t)+n](|u-y_2|-|u-y_1|)+[c_2(t)+n](|v-z_2|-|v-z_1|)\}\big|\\
\\
\leq& \sup_{(u,v)\in\R^{1+k}}\big|\{[c_1(t)+n]|u-y_2-u+y_1|+[c_2(t)+n]|v-z_2-v+z_1|\}\big|\\
\\
=& [c_1(t)+n]|y_1-y_2|+[c_2(t)+n]|z_1-z_2|.\nonumber\\
\end{split}
\end{equation}

(iv) For any $n\geq 1$, there exists $(u_n,v_n)\in \R^{1+k}$ such that
\begin{equation}
\begin{split}
f_n(t,y,z)&\leq f(t,u_n,v_n)-[c_1(t)+n]|u_n-y|-[c_2(t)+n]|v_n-z|+ n^{-1}.\nonumber
\end{split}
\end{equation}
In other words,
\[
f_n(t,y,z)+[c_1(t)+n]|u_n-y|+[c_2(t)+n]|v_n-z|\leq f(t,u_n,v_n)+ n^{-1}.
\]
Note that in order to make the left-hand side of above inequality finite as $n\rightarrow \infty$, it is necessary to have $\lim_{n\rightarrow \infty}(u_n,v_n)=(y,z)$. And then
\[
\lim_{n\rightarrow \infty}f_n(t,y,z)\leq f(t,y,z).
\]
On the other hand, by the definition of $f_n$, we have $f_n(t,y,z)\geq f(t,y,z)$. Hence
\[
\lim_{n\rightarrow \infty}f_n(t,y,z)=f(t,y,z).
\]\qed\\

Using functions $\{f_n\}_{n\geq1}$ as generators, we introduce the following sequence of equations
\begin{equation}
\label{bsde linear}
\widetilde{y}_n(t)= \xi+ \int_t^T f_n(s, \widetilde{y}_n(s),\widetilde{ z}_n(s))ds-\int_t^T \widetilde{z}_n(s)dW(s),\quad t\in[0,T].
\end{equation}

\begin{lemma} If conditions A3 hold, then equations \eqref{bsde linear} have unique solution pairs $(\widetilde{y}_n(\cdot),\widetilde{z}_n(\cdot))\in \widetilde{H}^2(0,T;\R)\times\widetilde{M}^2(0,T;\R^k)$, for any $n\geq1$.
\end{lemma}
\proof We only need to show that the assumptions of the previous section, which ensure the applicability of Theorem \ref{theorem1} (i), hold. Thus,
\begin{eqnarray}
\E\left[e^{\int_0^T\{\widetilde{\gamma(t)}+\widetilde{\beta_1}[c_1(t)+n]^2+\widetilde{\beta_2}[c_2(t)+n]^2\}dt}|\xi|^2\right]\leq e^{2Tn^2(\widetilde{\beta_1}+\widetilde{\beta_2})}\E[\widetilde{p}(T)|\xi|^2]<\infty,\nonumber
\end{eqnarray}
and
\begin{eqnarray}
&&\E\int_0^Te^{\int_0^t\{\widetilde{\gamma(t)}+\widetilde{\beta_1}[c_1(t)+n]^2+\widetilde{\beta_2}[c_2(t)+n]^2\}ds}\frac{|f_n(t,0,0)|^2}{\widetilde{\gamma(t)}+\widetilde{\beta_1}[c_1(t)+n]^2+\widetilde{\beta_2}[c_2(t)+n]^2}dt\nonumber\\
\nonumber\\
&&\leq\E\int_0^Te^{2n^2t(\widetilde{\beta_1}+\widetilde{\beta_2})}\widetilde{p}(t)\frac{|c_0(t)|^2}{\widetilde{\gamma(t)}+\widetilde{\beta_1}[c_1(t)+n]^2+\widetilde{\beta_2}[c_2(t)+n]^2}dt\nonumber\\
\nonumber\\
&&\leq2e^{2n^2T(\widetilde{\beta_1}+\widetilde{\beta_2})}\E\int_0^T\widetilde{p}(t)\frac{|c_0(t)|^2}{\widetilde{\alpha}(t)}dt<\infty.\nonumber
\end{eqnarray}
\qed\\

Our main task now is to prove that the sequence of solutions $\{\widetilde{y}_n(\cdot),\widetilde{z}_n(\cdot)\}_{n\geq 1}$, converges to the solution $(y(\cdot),z(\cdot))$ of \eqref{bsde}. We first present two useful lemmas.
\begin{lemma}
\label{bound}
Let the conditions A3 hold. There exists a constant $\kappa$, independent of $n$, such that $\widetilde{\|\widetilde{y}_n\|}\leq \kappa$ and $\widetilde{\|\widetilde{z}_n\|}\leq \kappa$, for all $n\geq1$.
\end{lemma}
\proof
From the fact that the sequence $f_n$ is decreasing, and Theorem \ref{compare} (i), we know that $\widetilde{y}_1(t)\geq\widetilde{y}_2(t)\geq...,$ $\forall t\in[0,T]$ $a.s.$. Hence, there exists a constant $\kappa_1$ such that
\begin{eqnarray}
\kappa_1\geq\widetilde{\|\widetilde{y}_1\|}\geq\widetilde{\|\widetilde{y}_2\|}\geq....\nonumber
\end{eqnarray}
By making use the linear growth property of $f_n$, we obtain
\begin{equation}
\begin{split}
&-d\widetilde{p}(t)|\widetilde{y}_n(t)|^2\\
\\
=&-\widetilde{\alpha}(t)\widetilde{p}(t)|\widetilde{y}_n(t)|^2dt-\widetilde{p}(t)|\widetilde{z}_n(t)|^2dt- 2\widetilde{p}(t)\widetilde{y}_n(t)\widetilde{z}_n(t)dW(t)\\
\\
&+2\widetilde{p}(t)\widetilde{y}_n(t)f_n(t,\widetilde{y}_n(t),\widetilde{z}_n(t))dt\\
\\
\leq&-\widetilde{\alpha}(t)\widetilde{p}(t)|\widetilde{y}_n(t)|^2dt-\widetilde{p}(t)|\widetilde{z}_n(t)|^2dt- 2\widetilde{p}(t)\widetilde{y}_n(t)\widetilde{z}_n(t)dW(t)\\
\\
&+2\widetilde{p}(t)|\widetilde{y}_n(t)||f_n(t,\widetilde{y}_n(t),\widetilde{z}_n(t))|dt\\
\\
\leq&-\widetilde{\alpha}(t)\widetilde{p}(t)|\widetilde{y}_n(t)|^2dt-\widetilde{p}(t)|\widetilde{z}_n(t)|^2dt- 2\widetilde{p}(t)\widetilde{y}_n(t)\widetilde{z}_n(t)dW(t)\\
\\
&+2\widetilde{p}(t)|\widetilde{y}_n(t)|[c_0(t)+c_1(t)|\widetilde{y}_n(t)|+c_2(t)|\widetilde{z}_n(t)|]dt\\
\\
\leq&-\widetilde{\alpha}(t)\widetilde{p}(t)|\widetilde{y}_n(t)|^2dt-\widetilde{p}(t)|\widetilde{z}_n(t)|^2dt- 2\widetilde{p}(t)\widetilde{y}_n(t)\widetilde{z}_n(t)dW(t)\\
\\
&+\widetilde{p}(t)|\widetilde{y}_n(t)|^2dt+\widetilde{p}(t)c^2_0(t)dt+\widetilde{\beta_1}c_1^2(t)\widetilde{p}(t)|\widetilde{y}_n(t)|^2dt+(\widetilde{\beta_1})^{-1}\widetilde{p}(t)|\widetilde{y}_n(t)|^2dt\\
\\
&+\widetilde{\beta_2}c_2^2(t)\widetilde{p}(t)|\widetilde{y}_n(t)|^2+(\widetilde{\beta_2})^{-1}\widetilde{p}(t)|\widetilde{z}_n(t)|^2dt\\
\\
\leq&-\left[1-(\widetilde{\beta_2})^{-1}\right]\widetilde{p}(t)|\widetilde{z}_n(t)|^2dt+\widetilde{p}(t)c^2_0(t)dt\\
\\
&+\left[1+(\widetilde{\beta_1})^{-1}\right]\widetilde{p}(t)|\widetilde{y}_n(t)|^2dt-2\widetilde{p}(t)\widetilde{y}_n(t)\widetilde{z}_n(t)dW(t),\nonumber
\end{split}
\end{equation}
which in integral form becomes
\begin{eqnarray}
\int_t^T\widetilde{p}(s)|\widetilde{z}_n(s)|^2ds&\leq& \frac{\widetilde{p}(T)\xi^2}{\left[1-(\widetilde{\beta_2})^{-1}\right]}+\frac{\int_t^T\widetilde{p}(s)c^2_0(s)ds}{\left[1-(\widetilde{\beta_2})^{-1}\right]}\nonumber\\
&+&\frac{\left[1+(\widetilde{\beta_1})^{-1}\right]}{\left[1-(\widetilde{\beta_2})^{-1}\right]}\int_t^T\widetilde{p}(s)|\widetilde{y}_n(t)|^2ds-\frac{2\int_t^T\widetilde{p}(s)\widetilde{y}_n(s)\widetilde{z}_n(s)dW(t)}{\left[1-(\widetilde{\beta_2})^{-1}\right]}.\nonumber
\end{eqnarray}
From Lemma \ref{first}, we know that the stochastic integral on the right-hand side is a martingale. Taking the expectation of both sides gives
\begin{eqnarray}
\widetilde{\|\widetilde{z}_n\|}&\leq& \frac{\E[\widetilde{p}(T)\xi^2]}{\left[1-(\widetilde{\beta_2})^{-1}\right]}+\frac{\E\int_0^T\widetilde{p}(s)c^2_0(s)ds}{\left[1-(\widetilde{\beta_2})^{-1}\right]}
+\frac{\left[1+(\widetilde{\beta_1})^{-1}\right]}{\left[1-(\widetilde{\beta_2})^{-1}\right]}\widetilde{\|\widetilde{y}_n\|}\\
&\leq& \frac{\E[\widetilde{p}(T)\xi^2]}{\left[1-(\widetilde{\beta_2})^{-1}\right]}+\frac{\E\int_0^T\widetilde{p}(s)c^2_0(s)ds}{\left[1-(\widetilde{\beta_2})^{-1}\right]}
+\frac{\left[1+(\widetilde{\beta_1})^{-1}\right]}{\left[1-(\widetilde{\beta_2})^{-1}\right]}\kappa_1=\kappa_2.\nonumber
\end{eqnarray}
Finally, $\kappa=\max(\kappa_1,\kappa_2)$.
\qed
\\
\begin{lemma}
\label{conv}
Let the conditions A3 hold. The pair of processes $(\widetilde{y}_n(\cdot),\widetilde{z}_n(\cdot))_{n\geq 1}$ converges to $(\widetilde{y}(\cdot),\widetilde{z}(\cdot))$ in ${\widetilde{M}}^2(0,T; \R)\times {\widetilde{M}}^2(0,T; \R^k)$.

\proof
Let us consider a measurable and Lipchitz function $g(t,y,z)=-[c_0(t)+c_1(t)|y|+c_2(t)|z|]$. From Lemma \ref{seq} (iii) and Theorem \ref{theorem1} (i) we know the following BSDEs have a unique adapted solution on ${\widetilde{H}}^2(0,T; \R)\times {\widetilde{M}}^2(0,T;\R^k)$:
\[
\widetilde{y}_n(t)= \xi+ \int_t^T f_n(s,\widetilde{y}_n(s),\widetilde{z}_n(s))ds-\int_t^T \widetilde{z}_n(s)dW(s),\quad t\in[0,T],
\]
and
\[
K(t)= \xi+ \int_t^T g(s,K(s),L(s))ds-\int_t^T L(s)dW(s),\quad t\in[0,T].
\]
By Comparison Theorem \ref{compare} (i), we have
\[
K(t)\leq \widetilde{y}_{n}(t)\leq \widetilde{y}_{n-1}(t)\leq \widetilde{y}_1(t),\,\, \forall n\geq 1.
\]
Hence $\{\widetilde{y}_n(t)\}_{n\geq 1}$ is decreasing and bounded in ${\widetilde{M}}^2(0,T; \R)$. Then by dominated convergence theorem, we know that $\{\widetilde{y}_n(t)\}_{n\geq 1}$ converges pointwisely to $\widetilde{y}(t)^*$ in ${\widetilde{M}}^2(0,T;\R)$.

Applying the It\^o's formula to $p(t)|{\widetilde{y}}_{n}(t)-{\widetilde{y}}_m(t)|^2$ and writing in integral form, we have
\begin{equation}
\label{mn}
\begin{split}
&\widetilde{p}(t)|\widetilde{y}_n(t)-\widetilde{y}_m(t)|^2+\int^T_t\widetilde{p}(s)|\widetilde{z}_n(s)-\widetilde{z}_m(s)|^2ds\\
\\
=&\int^T_t2\widetilde{p}(s)(\widetilde{y}_n(s)-\widetilde{y}_m(s))[f_n(s,\widetilde{y}_n(s),\widetilde{z}_n(s))-f_m(s,\widetilde{y}_m(s),\widetilde{z}_m(s))]ds\\
\\
+&\int_t^T-2\widetilde{p}(s)(\widetilde{y}_n(s)-\widetilde{y}_m(s))(\widetilde{z}_n(s)-\widetilde{z}_m(s))dW(s).\\
\end{split}
\end{equation}
Taking expectations on both sides, from Lemma \ref{first} (ii) and using the linear growth property of $f_n$ and $f_m$, we have
\begin{equation}
\begin{split}
&\E\left[\int^T_t\widetilde{p}(s)|\widetilde{z}_n(s)-\widetilde{z}_m(s)|^2ds\right]\\
\\
\leq&\,2\E\left[\int^T_t\widetilde{p}(s)(\widetilde{y}_n(s)-\widetilde{y}_m(s))[f_n(s,\widetilde{y}_n(s),\widetilde{z}_n(s))-f_m(s,\widetilde{y}_m(s),\widetilde{z}_m(s))]ds\right]\\
\\
\leq&\,2\E\left[\int^T_t\sqrt{\widetilde{\alpha}(s)\widetilde{p}(s)}(\widetilde{y}_n(s)-\widetilde{y}_m(s))\frac{\sqrt{\widetilde{p}(s)}}{\sqrt{\widetilde{\alpha}(s)}}[f_n(s,\widetilde{y}_n(s),\widetilde{z}_n(s))-f_m(s,\widetilde{y}_m(s),\widetilde{z}_m(s))]ds\right]\\
\\
\leq&\,2\left(\E\left[\int^T_t\widetilde{\alpha}(s)\widetilde{p}(s)|\widetilde{y}_n(s)-\widetilde{y}_m(s)|^2ds\right]\right)^{\frac{1}{2}}\\
\\
&\cdot\left(\E\left[\int^T_t\frac{\widetilde{p}(s)}{\widetilde{\alpha}(s)}|f_n(s,\widetilde{y}_n(s),\widetilde{z}_n(s))-f_m(s,\widetilde{y}_m(s),\widetilde{z}_m(s))|^2ds\right]\right)^{\frac{1}{2}}\\
\\
\leq&\,2\sqrt{2}\left(\E\left[\int^T_t\widetilde{\alpha}(s)\widetilde{p}(s)|\widetilde{y}_n(s)-\widetilde{y}_m(s)|^2ds\right]\right)^{\frac{1}{2}}\bigg(\E\bigg[\int^T_t\frac{\widetilde{p}(s)}{\widetilde{\alpha}(s)}\bigg[|c_0(s)+c_1(s)|\widetilde{y}_n(s)|\\
\\
&+c_2(s)|\widetilde{z}_n(s)||^2+|c_0(s)+c_1(s)|\widetilde{y}_m(s)|+c_2(s)|\widetilde{z}_m(s)||^2\bigg]ds\bigg]\bigg)^{\frac{1}{2}}\\
\\
\leq&\,2\sqrt{6}\left(\E\left[\int^T_t\widetilde{\alpha}(s)\widetilde{p}(s)|\widetilde{y}_n(s)-\widetilde{y}_m(s)|^2ds\right]\right)^{\frac{1}{2}}\bigg(\E\bigg[\int^T_t\frac{\widetilde{p}(s)}{\widetilde{\alpha}(s)}\bigg[2c_0^2(s)+c_1^2(s)|\widetilde{y}_n(s)|^2\\
\\
&+c_2^2(s)|\widetilde{z}_n(s)|^2+c_1^2(s)|\widetilde{y}_m(s)|^2+c_2^2(s)|\widetilde{z}_m(s)|^2\bigg]ds\bigg]\bigg)^{\frac{1}{2}}\nonumber\\
\end{split}
\end{equation}
\begin{equation}
\begin{split}
\leq&\,2\sqrt{6}\left(\E\left[\int^T_t\widetilde{\alpha}(s)\widetilde{p}(s)|\widetilde{y}_n(s)-\widetilde{y}_m(s)|^2ds\right]\right)^{\frac{1}{2}}\bigg(\E\bigg[\int^T_t\widetilde{p}(s)\bigg[\frac{2c_0^2(s)}{\widetilde{\alpha}(s)}+|\widetilde{y}_n(s)|^2\\
\\
&+|\widetilde{z}_n(s)|^2+|\widetilde{y}_m(s)|^2+|\widetilde{z}_m(s)|^2\bigg]ds\bigg]\bigg)^{\frac{1}{2}}\\
\\
\leq&\,2\sqrt{6}\,\widetilde{\kappa}\left(\E\left[\int^T_t\widetilde{\alpha}(s)\widetilde{p}(s)|\widetilde{y}_n(s)-\widetilde{y}_m(s)|^2ds\right]\right)^{\frac{1}{2}},\nonumber\\
\end{split}
\end{equation}
where $\widetilde{\kappa}\equiv \left(4\kappa+ {\widetilde{\left\|c_0(\cdot){\alpha}(\cdot)^{-\frac{1}{2}}\right\|}}\right)^{\frac{1}{2}}$. Therefore, this, together with the fact that $\{\widetilde{y}_n(t)\}_{n\geq 1}$ pointwisely converges in ${\widetilde{M}}^2(0,T; \R)$, implies that $\{{\widetilde{z}}_n(t)\}_{n\geq 1}$ is a Cauchy sequence in
${\widetilde{M}}^2(0,T;\R^k)$ and then converges to ${\widetilde{z}}(t)^*$ in the same space.
\qed
\end{lemma}

Now we present the main result in this section.
\begin{theorem}
\label{cont} (Existence)
Equation \eqref{bsde} has an adapted solution $(\widetilde{y}(\cdot),\widetilde{z}(\cdot))\in {\widetilde{M}}^2(0,T; \R)\times {\widetilde{M}}^2(0,T; \R^k)$, which is also a maximal solution, i.e. for any other solution $(\bar{y}(\cdot),\bar{z}(\cdot))$ of equation \eqref{bsde}, we have $\widetilde{y}(\cdot)\geq \bar{y}(\cdot)$.
\end{theorem}
\proof
Similar to previous calculation in Lemma \ref{conv}, taking supremum over $t$ for equation \eqref{mn} and using the Burkholder-Davis-Gundy's inequality, we have
\begin{equation}
\begin{split}
&\E\left[\sup_{t\in[0,T]}\widetilde{p}(t)|\widetilde{y}_n(t)-\widetilde{y}_m(t)|^2\right]\\
\\
\leq&\,\E\left[\sup_{t\in[0,T]}\int^T_t2\widetilde{p}(s)(\widetilde{y}_n(s)-\widetilde{y}_m(s))[f_n(s,\widetilde{y}_n(s),\widetilde{z}_n(s))-f_m(s,\widetilde{y}_m(s),\widetilde{z}_m(s))]ds\right]\\
\\
&+\,\E\left[\sup_{t\in[0,T]}\int_t^T-2\widetilde{p}(s)(\widetilde{y}_n(s)-\widetilde{y}_m(s))(\widetilde{z}_n(s)-\widetilde{z}_m(s))dW(s)\right]\\
\\
\leq&\,2\E\left[\int^T_0\widetilde{p}(s)|\widetilde{y}_n(s)-\widetilde{y}_m(s)||f_n(s,\widetilde{y}_n(s),\widetilde{z}_n(s))-f_m(s,\widetilde{y}_m(s),\widetilde{z}_m(s))|ds\right]\\
\\
&+\,K\E\left[\int_0^T|\sqrt{\widetilde{p}(s)}(\widetilde{y}_n(s)-\widetilde{y}_m(s))|^2|\sqrt{\widetilde{p}(s)}(\widetilde{z}_n(s)-\widetilde{z}_m(s))|^2ds\right]^{\frac{1}{2}}\\
\\
\leq&\,2\left(\E\left[\int^T_t\widetilde{\alpha}(s)\widetilde{p}(s)|\widetilde{y}_n(s)-\widetilde{y}_m(s)|^2ds\right]\right)^{\frac{1}{2}}\\
\\
&\cdot\left(\E\left[\int^T_t\frac{\widetilde{p}(s)}{\widetilde{\alpha}(s)}|f_n(s,\widetilde{y}_n(s),\widetilde{z}_n(s))-f_m(s,\widetilde{y}_m(s),\widetilde{z}_m(s))|^2ds\right]\right)^{\frac{1}{2}}\\
\\
&+\,\frac{K}{2}\E\left[\sup_{t\in [0,T]}|\sqrt{\widetilde{p}(s)}(\widetilde{y}_n(s)-\widetilde{y}_m(s))|^2+\int_t^T|\sqrt{\widetilde{p}(s)}(\widetilde{z}_n(s)-\widetilde{z}_m(s))|^2ds\right]\\
\\
\leq&\,2\sqrt{6}\,\widetilde{\kappa}\left(\E\left[\int^T_t\widetilde{\alpha}(s)\widetilde{p}(s)|\widetilde{y}_n(s)-\widetilde{y}_m(s)|^2ds\right]\right)^{\frac{1}{2}}\\
\\
&+\,\frac{K}{2}\E\left[\sup_{t\in [0,T]}\widetilde{p}(s)|\widetilde{y}_n(s)-\widetilde{y}_m(s)|^2+\int_t^T\widetilde{p}(s)|\widetilde{z}_n(s)-\widetilde{z}_m(s)|^2ds\right].\nonumber\\
\end{split}
\end{equation}
Therefore by Lemma \ref{conv}, for any $t\in[0,T]$, $\{\widetilde{y}_n(t)\}_{n\geq 1}$ converges uniformly in ${\widetilde{M}}^2(0,T;\R)$ to $y(t)$ in the same space. Since $\{\widetilde{y}_n(t)\}_{n\geq 1}$ is continuous, by the uniform convergence theorem, $\widetilde{y}(t)$ is a continuous process.

Now we show that the sequence of processes $\{f_n(t,\widetilde{y}_n(t),\widetilde{z}_{n}(t))\}_{n\geq 1}$ converges to $\{f(t,\widetilde{y}(t),\widetilde{z}(t))\}$ in
$L_{\mathscr{F}}(0,T;\R)$. Note that for any positive and large enough $\delta$,
\begin{eqnarray}
&& \E\left[\int_0^T|f_n(s,\widetilde{y}_n(s),\widetilde{z}_{n}(s))-f(s,\widetilde{y}(s),\widetilde{z}(s))|ds\right]\nonumber\\
\nonumber\\
&=& \E\left[\int_0^T|f_n(s,\widetilde{y}_n(s),\widetilde{z}_{n}(s))-f(s,\widetilde{y}(s),\widetilde{z}(s))|\mathds{1}_{\left\{\frac{c_1(s)|\widetilde{y}_n(s)|+c_2(s)|\widetilde{z}_n(s)|}{c_0^2(s)+c_1^2(s)+c_2^2(s)}\leq\delta\right\}}ds\right]\nonumber\\
\nonumber\\
&& +\:\E\left[\int_0^T|f_n(s,\widetilde{y}_n(s),\widetilde{z}_{n}(s))-f(s,\widetilde{y}(s),\widetilde{z}(s))|\mathds{1}_{\left\{\frac{c_1(s)|\widetilde{y}_n(s)|+c_2(s)|\widetilde{z}_n(s)|}{c_0^2(s)+c_1^2(s)+c_2^2(s)}>\delta\right\}}ds\right]\nonumber\\
\nonumber\\
&\leq& \E\left[\int_0^T|f_n(s,\widetilde{y}_n(s),\widetilde{z}_{n}(s))-f(s,\widetilde{y}_n(s),\widetilde{z}_n(s))|\mathds{1}_{\left\{\frac{c_1(s)|\widetilde{y}_n(s)|+c_2(s)|\widetilde{z}_n(s)|}{c_0^2(s)+c_1^2(s)+c_2^2(s)}\leq\delta\right\}}ds\right]\nonumber\\
\nonumber\\
&& +\E\left[\int_0^T|f(s,\widetilde{y}_n(s),\widetilde{z}_{n}(s))-f(s,\widetilde{y}(s),\widetilde{z}(s))|\mathds{1}_{\left\{\frac{c_1(s)|\widetilde{y}_n(s)|+c_2(s)|\widetilde{z}_n(s)|}{c_0^2(s)+c_1^2(s)+c_2^2(s)}\leq\delta\right\}}ds\right]\nonumber\\
\nonumber\\
&& +\E\left[\int_0^T|f_n(s,\widetilde{y}_n(s),\widetilde{z}_{n}(s))-f(s,\widetilde{y}(s),\widetilde{z}(s))|\mathds{1}_{\left\{\frac{c_1(s)|\widetilde{y}_n(s)|+c_2(s)|\widetilde{z}_n(s)|}{c_0^2(s)+c_1^2(s)+c_2^2(s)}>\delta\right\}}ds\right].\nonumber\\
\end{eqnarray}
By (ii), (iii) and (iv) in Lemma \ref{seq}, assumption (i)' and the Dini's Theorem, as
$n\rightarrow\infty$, we have
\[
\sup_{\left\{\frac{c_1(s)|\widetilde{y}_n(s)|+c_2(s)|\widetilde{z}_n(s)|}{c_0^2(s)+c_1^2(s)+c_2^2(s)}\leq\delta\right\}}|f_n(t,\widetilde{y}(t),\widetilde{z}(t))-f(t,\widetilde{y}(t),\widetilde{z}(t))|\longrightarrow 0.
\]
Therefore by the dominated convergence theorem, the first term in right hand side uniformly converges to $0$. Due to assumption (i)', at least along a subsequence, the second term in the right hand side converges to $0$. For the final term in the right hand side, by Lemma \ref{seq} (iii), assumption (ii)' and Lemma \ref{bound}, together with the fact that $a\mathds{1}_{\{X>a\}}<X$ for any nonnegative random variable $X$ and $a>0$, we have
\begin{equation}
\begin{split}
&\E\left[\int_0^T|f_n(t,\widetilde{y}_n(s),\widetilde{z}_{n}(s))-f(t,\widetilde{y}(s),\widetilde{z}(s))|\mathds{1}_{\left\{\frac{c_1(s)|\widetilde{y}_n(s)|+c_2(s)|\widetilde{z}_n(s)|}{c_0^2(s)+c_1^2(s)+c_2^2(s)}>\delta\right\}}ds\right]\\
\\
\leq&\,\E\left[\int_0^T\left[2c_0(s)+c_1(s)|\widetilde{y}_n(s)|+c_2(s)|\widetilde{z}_n(s)|+c_1(s)|\widetilde{y}(s)|+c_2(s)|\widetilde{z}(s)\right]
\mathds{1}_{\left\{\frac{c_1(s)|\widetilde{y}_n(s)|+c_2(s)|\widetilde{z}_n(s)|}{c_0^2(s)+c_1^2(s)+c_2^2(s)}>\delta\right\}}ds\right]\\
\\
\leq&\,\E\bigg[\int_0^T[2c_0(s)+c_1(s)|\widetilde{y}_n(s)|+c_2(s)|\widetilde{z}_n(s)|+c_1(s)|\widetilde{y}(s)|+c_2(s)|\widetilde{z}(s)]\\
&\quad\quad\frac{c_1(s)|\widetilde{y}_n(s)|+c_2(s)|\widetilde{z}_n(s)|}{\delta(c_0^2(s)+c_1^2(s)+c_2^2(s))}ds\bigg]\\
\\
=& \, \frac{1}{\delta}\,\E\bigg[\int_0^T\frac{1}{c_0^2(s)+c_1^2(s)+c_2^2(s)}[2c_0(s)c_1(s)|\widetilde{y}_n(s)|+c_1^2(s)|\widetilde{y}_n(s)|^2+c_1(s)c_2(s)|\widetilde{y}_n(s)||\widetilde{z}_n(s)|\\
\\
&+c_1^2(s)|\widetilde{y}_n(s)||\widetilde{y}(s)|+c_1(s)c_2(s)|\widetilde{y}_n(s)||\widetilde{z}(s)|+2c_0(s)c_2(s)|\widetilde{z}_n(s)|+c_1(s)c_2(s)|\widetilde{y}_n(s)||\widetilde{z}_n(s)|\\
\\
&+c_2^2(s)|\widetilde{z}_n(s)|^2+c_1(s)c_2(s)|\widetilde{y}(s)||\widetilde{z}_n(s)|+c_2^2(s)|\widetilde{z}_n(s)||\widetilde{z}(s)|]ds\bigg]\\
\\
\leq& \,\frac{1}{\delta}\,\E\bigg[\int_0^T\frac{1}{c_0^2(s)+c_1^2(s)+c_2^2(s)}[2c_0^2(s)+4c_1^2(s)|\widetilde{y}_n(s)|^2+4c_2^2(s)|\widetilde{z}_n(s)|^2\\
\\
&+c_1^2(s)|\widetilde{y}(s)|^2+c_2^2(s)|\widetilde{z}(s)|^2]ds\bigg]\\
\\
\leq&\,\frac{1}{\delta}\,\E\bigg[\int_0^T[2+4|\widetilde{y}_n(s)|^2+4|\widetilde{z}_n(s)|^2+|\widetilde{y}(s)|^2+|\widetilde{z}(s)|^2]ds\bigg]\leq \frac{\bar{\kappa}}{\delta},\nonumber\\
\\
\end{split}
\end{equation}
where $\bar{\kappa}$ is a constant independent of $n$.
Hence taking limits in the following equation
\begin{equation}
\widetilde{y}_n(t)= \xi+ \int_t^T f_n(s,\widetilde{y}_n(s),\widetilde{z}_n(s))ds-\int_t^T \widetilde{z}_n(s)dW(s),\quad t\in[0,T]\nonumber.
\end{equation}
we deduce that $(\widetilde{y}(t),\widetilde{z}(t))\in {\widetilde{M}}^2(0,T;\R)\times {\widetilde{M}}^2(0,T;\R^k)$ is an adapted solution of equation \eqref{bsde}.

Furthermore, suppose that $(\bar{y}(t),\bar{z}(t))\in {\widetilde{M}}^2(0,T;\R)\times {\widetilde{M}}^2(0,T;\R^k)$ is any solution of equation \eqref{bsde}. By the comparison theorem
\ref{compare}, we have $\widetilde{y}_n(t)\geq \bar{y}(t)$ for any $n\geq 1$ and then $\widetilde{y}(t)\geq \bar{y}(t)$. Hence $\widetilde{y}(t)$ is a maximal solution of equation
\eqref{bsde}.
\qed
\\

\section{Conclusions} We have considered two classes of BSDEs with possibly unbounded generator. The first class has a Lipschitz-type generator,
whereas the second class has a continuous generator that satisfies a certain linear growth condition. In both cases we have given sufficient conditions
for solvability. We expect that these results will prove to be useful in tackling more difficult problems with unbounded generator, such as the BSDEs
with a quadratic growth and the Riccati BSDE, which play a fundamental role in stochastic control.

\section*{Appendix}
Here we include the derivation of the lower bound for the parameter $\beta$ that appears in~\cite{KH}. We do so for the readers' convenience, since in Theorem 6.1 of~\cite{KH} no explicit lower bound is given, it is only assumed that parameter $\beta$ should be {\it large enough}. The notation of~\cite{KH} will be used.\\

Equation (6.5) of~\cite{KH} states that for some constants $k$ and $k'$ the following holds
\begin{eqnarray}
\|(y,\eta)\|^2_{\beta}=k\|\xi\|^2_{\beta}+\frac{k'}{\beta}\left\|\frac{f}{\alpha}\right\|^2_{\beta},\label{1}
\end{eqnarray}
where the definitions of these norms are given in~\cite{KH}, and are just weighted Euclidian norms. From equation (5.5) of~\cite{KH}, which gives the definition of the norm $\|(y,\eta)\|^2_{\beta}$, and the conclusions of Lemma 6.2 of~\cite{KH}, we obtain
\begin{eqnarray}
\|(y,\eta)\|^2_{\beta}&=&\, \|\alpha y\|^2_{\beta}+ \|\eta\|^2_{\beta}\nonumber\\
&\leq&\, \frac{2}{\beta}\|\xi\|^2_{\beta} + \frac{8}{\beta^2}\left\|\frac{f}{\alpha}\right\|^2_{\beta}+ 18\|\xi\|^2_{\beta}+\frac{45}{\beta}\left\|\frac{f}{\alpha}\right\|^2_{\beta}\nonumber\\
&=& \left(18+\frac{2}{\beta}\right)\|\xi\|^2_{\beta}+\left(\frac{45}{\beta}+\frac{8}{\beta^2}\right)\left\|\frac{f}{\alpha}\right\|^2_{\beta}.\label{2}
\end{eqnarray}
Comparing (\ref{1}) and (\ref{2}) gives $k'= 45+\frac{8}{\beta}$.\\

Equation (6.16) of~\cite{KH} states that for some constants $\tilde{k}$ and $\tilde{k}'$ the following holds
\begin{eqnarray}
\|(\delta Y,\delta Z, \delta N)\|^2_{\beta}\leq&\, \tilde{k}\|\delta \xi\|^2_{\beta}+ \frac{\tilde{k}'}{\beta}\left\|\frac{\delta_2f}{\alpha}\right\|^2_{\beta}.\label{3}
\end{eqnarray}
Here, different from~\cite{KH}, we have used the {\it tilde} notation for the constants $k$ and $k'$ in order to avoid the clash of notation with these constants introduced in the previous paragraph. By inequality (6.5) of~\cite{KH}, we obtain that
\begin{eqnarray}
\|(\delta Y,\delta Z, \delta N)\|^2_{\beta}&\leq&\, k\|\delta \xi\|^2_{\beta}+ \frac{k'}{\beta}\left\|\frac{\varphi_t}{\alpha^2_t}\right\|^2_{\beta}\nonumber\\
&\leq&\, k\|\delta \xi\|^2_{\beta}+ \frac{3k'}{\beta}\left(\left\|\alpha \delta Y\right\|^2_{\beta}+\left\|m^*\delta Z\right\|^2_{\beta}+\left\|\frac{\delta_2 f}{\alpha}\right\|^2_{\beta}\right)\nonumber\\
&\leq&\, k\|\delta \xi\|^2_{\beta}+ \frac{3k'}{\beta}\left(k\left\|\delta \xi\right\|^2_{\beta}+\frac{k'}{\beta}\left\|\frac{\delta_2 f}{\alpha}\right\|^2_{\beta}+\left\|\frac{\delta_2 f}{\alpha}\right\|^2_{\beta}\right)\nonumber\\
&=&\, \left(k+\frac{3kk'}{\beta}\right)\|\delta \xi\|^2_{\beta}+ \frac{3k'}{\beta}\left(\frac{k'}{\beta}+1\right)\left\|\frac{\delta_2 f}{\alpha}\right\|^2_{\beta}.\label{4}
\end{eqnarray}
Comparing (\ref{3}) and (\ref{4}) gives $\widetilde{k'}= 3k'\left(\frac{k'}{\beta}+1\right)$.\\

The inequality at the end of page 35 of~\cite{KH} is
\begin{eqnarray}
\|\alpha \delta Y\|^2_{\beta}+ \|m^*\delta Z\|^2_{\beta}\leq \frac{\widehat{k'}}{\beta}\|\alpha \delta y\|^2_{\beta}+ \|m^*\delta z\|^2_{\beta},\label{5}
\end{eqnarray}
where we have used the {\it hat} notation for the constant $k'$ in order to avoid the clash of notation with this constant introduced earlier.
Similarly to the previous paragraph we obtain
\begin{eqnarray}
\|(\delta Y, \delta Z)\|^2_{\beta}&=&\, \|\alpha \delta Y\|^2_{\beta}+ \|m^*\delta Z\|^2_{\beta}\leq\, \frac{\widetilde{k'}}{\beta}\left\|\frac{\varphi}{\alpha}\right\|^2_{\beta}\nonumber\\
&\leq&\, \frac{3\widetilde{k'}}{\beta}\|\alpha \delta y\|^2_{\beta}+ \|m^*\delta z\|^2_{\beta}.\label{6}
\end{eqnarray}
Comparing (\ref{5}) and (\ref{6}) gives $\widehat{k'}=3\widetilde{k}'=9k'\left(\frac{k'}{\beta}+1\right)$. In order to apply the contraction mapping principle, it is necessary to have $\frac{\widehat{k'}}{\beta}<1$, i.e.
\[
9\left(\frac{45}{\beta}+\frac{8}{\beta^2}\right)\left(\frac{45}{\beta}+\frac{8}{\beta^2}+1\right)<1.
\]
By solving above inequality for $\beta>0$, we obtain that by {\it large enough} in~\cite{KH} it is meant that $\beta> 446.05$.\\

\section*{Acknowledgement}
We are grateful to the reviewer for a careful reading of the paper and suggestions that led to an improved version.

\end{document}